\documentclass[11pt]{article}

\usepackage{a4,amssymb,amsmath,color,graphicx}
\usepackage{trfsigns}

\def\eins{\mbox{1\hskip-0.24em l}}
\def\T{^{\sf T}}
\def\mT{^{\sf -T}}
\newcommand{\C}{ {\mathbb C} }

\newcommand{\R}{ {\mathbb R} }

\newcommand{\MM}{{\mathbb M}}
\newcommand{\Uu}{{\mathbb U}}
\newcommand{\diag}{\,\mbox{diag}}

\newcommand{\cc}{{\bf c}}
\newcommand{\LL}{{\cal L}}
\newcommand{\PP}{{\cal P}}
\newcommand{\gdw}{\ \iff\ }

\newcommand{\qed}{\qquad\mbox{$\square$}}

\newtheorem{remark}{Remark}[section]
\newtheorem{theorem}{Theorem}[section]
\newtheorem{lemma}{Lemma}[section]

\begin{document}
\title{Implicit A-Stable Peer Triplets for ODE Constrained Optimal Control
Problems}
\author{Jens Lang \\
{\small \it Technical University Darmstadt,
Department of Mathematics} \\
{\small \it Dolivostra{\ss}e 15, 64293 Darmstadt, Germany}\\
{\small lang@mathematik.tu-darmstadt.de} \\ \\
Bernhard A. Schmitt \\
{\small \it Philipps-Universit\"at Marburg,
Department of Mathematics,}\\
{\small \it Hans-Meerwein-Stra{\ss}e 6, 35043 Marburg, Germany} \\
{\small schmitt@mathematik.uni-marburg.de}}
\maketitle

\begin{abstract}
This paper is concerned with the construction and convergence analysis of
novel implicit Peer triplets of two-step nature with four stages for nonlinear
ODE constrained optimal control problems. We combine the property of superconvergence
of some standard Peer method for inner grid points with carefully designed starting and
end methods to achieve order four for the
state variables and order three for the adjoint variables in a
first-discretize-then-optimize approach together with A-stability.
The notion triplets emphasizes that these three different Peer methods have to
satisfy additional matching conditions.
Four such Peer triplets of practical interest are constructed.
Also as a benchmark method, the well-known backward differentiation formula BDF4, which is only
$A(73.35^o)$-stable,
is extended to a special Peer triplet to supply an adjoint consistent method of higher order
and BDF type with equidistant nodes. Within the class of Peer triplets, we found
a diagonally implicit $A(84^o)$-stable method with nodes symmetric in $[0,1]$ to a common
center that performs equally well. Numerical tests with three well
established optimal control problems confirm the theoretical findings also concerning A-stability.
\end{abstract}


\newpage

\section{Introduction}
The design of efficient time integrators for the numerical solution of optimal control
problems constrained by systems of ordinary differential equations (ODEs) is still an active
research field. Such systems typically arise from semi-discretized partial differential equations
describing, e.g., the dynamics of heat and mass transfer or fluid flow in complex physical systems.
Symplectic one-step Runge-Kutta methods \cite{LiuFrank2021,SanzSerna2016} exploit
the Hamiltonian structure of the first-order optimality system - the necessary conditions to find
an optimizer - and automatically yield a consistent approximation of the adjoint equations, which can
be used to compute the gradient of the objective function. The first-order symplectic Euler,
second-order St\"ormer-Verlet and higher-order Gauss methods are prominent representatives of this class,
which are all implicit for general Hamiltonian systems, see the monograph \cite{HairerWannerLubich2006}.
Generalized partitioned Runge-Kutta methods which allow to compute exact gradients with respect to the
initial condition are studied in \cite{MatsudaMiyatake2021}.
To avoid the solution of large systems of nonlinear equations, stabilized explicit Runge-Kutta-Chebyshev
methods have been recently proposed, too \cite{AlmuslimaniVilmart2021}. However, as all one-step methods,
also symplectic Runge-Kutta schemes join the structural suffering of order reductions, which
may lead to poor results in their application, e.g., to boundary control problems such as external cooling and heating in a
manufacturing process; see \cite{LubichOstermann1995,OstermannRoche1992} for a detailed study of this
behaviour.
\par
In contrast, multistep methods including Peer two-step methods avoid order reductions and allow a simple
implementation \cite{GerischLangPodhaiskyWeiner2009,GottermeierLang2009}. However, the discrete adjoint schemes
of linear multistep methods are in general not consistent or show a significant decrease of their approximation order \cite{AlbiHertyPareschi2019,Sandu2008}.
Recently, we have developed implicit Peer two-step methods \cite{LangSchmitt2020} with three stages
to solve ODE constrained optimal control problems of the form
\begin{align}
\mbox{minimize } C\big(y(T)\big) \label{OCprob_objfunc} &\\
\mbox{subject to } y'(t) =& \,f\big(y(t),u(t)\big),\quad
u(t)\in U_{ad},\;t\in(0,T], \label{OCprob_ODE}\\
y(0) =& \,y_0, \label{OCprob_ODEinit}
\end{align}
with the state $y(t)\in\R^m$, the control $u(t)\in\R^d$,
$f: \R^m\times\R^d\mapsto\R^m$, the objective function $C: \R^m\mapsto\R$, where the set of admissible controls $U_{ad}\subset\R^d$ is closed and convex. Introducing for any $u\in U_{ad}$ the normal cone mapping
\begin{align}
\label{def_cone}
N_U(u) =&\, \{ w\in\R^d: w^T(v-u)\le 0 \mbox{ for all } v\in U_{ad}\},
\end{align}
the first-order Karush-Kuhn-Tucker (KKT) optimality system \cite{Hager2000,Troutman1996} reads
\begin{align}
y'(t) =& \,f\big(y(t),u(t)\big),\quad t\in(0,T],\quad y(0)=y_0, \label{KKT_state}\\
p'(t) =& \,-\nabla_y f\big(y(t),u(t)\big)\T p(t),\quad t\in[0,T),
\quad p(T)=\nabla_y C\big(y(T)\big)\T, \label{KKT_adj}\\
& \,-\nabla_u f\big(y(t),u(t)\big)\T p(t) \in N_U\big(u(t)\big),\quad t\in[0,T]. \label{KKT_ctr}
\end{align}
In this paper, we assume the existence of a unique local solution $(y^\star,p^\star,u^\star)$
of the KKT system with sufficient
regularity properties to justify the use of higher order Peer triplets, see, e.g., the smoothness
assumption in \cite[Section 2]{Hager2000}.
\par
Following a \textit{first-discretize-and-then-optimize} approach, we apply an $s$-stage
implicit Peer two-step method to
\eqref{OCprob_ODE}-\eqref{OCprob_ODEinit} with approximations
$Y_{ni}\approx y(t_n+c_ih)$ and $U_{ni}\approx u(t_n+c_ih)$,
$i=1,\ldots,s,$ on an equi-spaced time grid $\{t_0,\ldots,t_{N+1}\}\subset [0,T]$ with
step size $h=(T-t_0)/(N+1)$ and nodes $c_1,\ldots,c_s$, which are fixed for all time steps, to
get the discrete constraint nonlinear optimal
control problem
\begin{align}
\mbox{minimize } C\big(y_h(T)\big) \label{OCprob_disc_objfunc} &\\[1mm]
\mbox{subject to } A_0 Y_0=&\,a\otimes y_0+hb\otimes f(y_0,u_0)+hK_0F(Y_0,U_0),\label{OCprob_peer_init}\\[1mm]
A_n Y_n=&\,B_nY_{n-1}+hK_nF(Y_n,U_n),\ n=1,\ldots,N,\label{OCprob_peer}\\[1mm]
y_h(T)=&\,(w^T\otimes I)Y_N,\label{OCprob_peer_yh}
\end{align}
with long vectors $Y_n=(Y_{ni})_{i=1}^s\in\R^{sm}$,
$U_n=(U_{ni})_{i=1}^s\in\R^{sd}$, and $F(Y_n,U_n)=\big(f(Y_{ni},U_{ni})\big)_{i=1}^s$.
Further, $y_h(T)\approx y(T)$, $u_0\approx u(0)$, and $a,b,w\in\R^s$ are additional
parameter vectors at both interval ends,
$A_n,B_n,K_n\in\R^{s\times s}$, and
$I\in\R^{m\times m}$ being the identity matrix. We will use the same symbol for a
coefficient matrix like $A$ and its Kronecker product $A\otimes I$ as a mapping from
the space $\R^{sm}$ to itself. In contrast to one-step methods, Peer two-step
methods compute $Y_n$ from the previous stage vector $Y_{n-1}$. Hence, also a
starting method, given in \eqref{OCprob_peer_init}, for the first time interval
$[t_0,t_1]$ is required. On each subinterval,
Peer methods may be defined by three coefficient matrices $A_n,B_n,K_n$, where $A_n$ and
$K_n$ are assumed to be nonsingular. For practical reasons, this general version will not
be used, the coefficients in the inner grid points will belong to a fixed Peer method
$(A_n,B_n,K_n)\equiv (A,B,K)$, $n=1,\ldots,N-1$. The last forward step has the same form
as the standard steps but uses exceptional coefficients $(A_N,B_N,K_N)$ to allow for a
better approximation of the end conditions.
\par
The KKT conditions \eqref{KKT_state}-\eqref{KKT_ctr} for ODE-constrained optimal control problems on a time interval
$[0,T]$ lead to a boundary value problem for a system of two differential equations, see
Section \ref{sec:bvp}.
The first one corresponds to the original \textit{forward} ODE for the state solution $y(t)$
and the second one is a linear, \textit{adjoint} ODE for a Lagrange multiplier $p(t)$.
It is well known that numerical methods for such problems may have to satisfy additional order conditions for the adjoint equation \cite{BonnansLaurentVarin2006,Hager2000,HertyPareschiSteffensen2013,LangVerwer2013,Sandu2006,SanzSerna2016}.
While these additional conditions are rather mild for one-step methods they may lead to severe restrictions for other types of methods like multistep and Peer methods, especially at the right-hand boundary at the end point $T$. Here, the order for the approximation of the adjoint equation may be limited to one.
\par
For Peer methods, this question was discussed first in \cite{SchroederLangWeiner2014} and the adjoint boundary condition at $T$ was identified as the most critical point.
In a more recent article \cite{LangSchmitt2020}, these bottlenecks could be circumvented by two measures.
First, equivalent formulations of the forward method are not equivalent for the adjoint formulation and using a redundant formulation of Peer methods with three coefficient matrices $(A,B,K)$ adds additional free parameters.
The second measure is to consider different methods for the first and last time interval.
Hence, instead of one single Peer method (which will be called the standard method) we discuss {\em triplets} of Peer methods consisting of a common standard method $(A,B,K)$ for all subintervals of the grid from the interior of $[0,T]$, plus a starting method $(A_0,K_0)$ for the first subinterval and an end method $(A_N,B_N,K_N)$ for the last one.
These two boundary methods may have lower order than the standard method since they are used only once.
\par
The present work extends the results from \cite{LangSchmitt2020} which considered methods with $s=3$ stages only, in two ways.
We will now concentrate on methods with four stages and better stability properties like A-stability.
The purpose of an accurate solution of the adjoint equation increases the number of conditions on the parameters of the method.
Requiring high order of convergence $s$ for the state variable $y(t)$ and order $s-1$ for the adjoint variable $p(t)$ -- which we combine to the pair $(s,s-1)$ from now on -- a variant of the method BDF3 was identified in \cite{LangSchmitt2020} as the most attractive standard method.
However, this method is not A-stable, with a stability angle of $\alpha=86^o$.
In order to obtain A-stability, we will reduce the required orders by one.
Still, we will show that convergence with the orders $(s,s-1)$ can be regained by a superconvergence property.
\par
The paper is organized as follows: In Section 2, the boundary value problem arising from the KKT system by
eliminating the control, and its discretization by means of discrete adjoint Peer two-step triplets are formulated.
An extensive error analysis concentrating on the superconvergence effect is presented in Section 3.
The restrictions imposed by the starting and end method on the standard Peer two-step method is studied in Section 4. The following Section 5 describes the actual construction principles of Peer triplets. Numerical tests are done in Section 6. The paper concludes with a summary in Section 7.

\section{The Boundary Value Problem}\label{sec:bvp}
Following the usual Lagrangian approach applied in \cite{LangSchmitt2020}, the first order discrete optimality conditions now consist of the forward equations \eqref{OCprob_peer_init}-\eqref{OCprob_peer_yh},
the discrete adjoint equations, acting backwards in time,
\begin{align}
A_N\T P_N=&\,w\otimes p_h(T)+h\nabla_YF(Y_N,U_N)\T K_N\T P_N,\label{KKT_adj_peer_init}\\[1mm]
A_n\T P_n=&\,B_{n+1}\T P_{n+1}+h\nabla_YF(Y_n,U_n)\T K_n\T P_n,\ N-1\ge n\ge 0,\label{KKT_adj_peer}
\end{align}
and the control conditions
\begin{align}
&\,-h\nabla_UF(Y_n,U_n)\T K_n\T P_n \in N_{U^s}\big(U_n\big),\ 0 \le n\le N,\label{KKT_ctr_peer}\\[1mm]
&\,-h\nabla_{u_0}f(y_0,u_0)\T (b\T\otimes I) P_0 \in N_U\big(u_0\big).\label{KKT_ctr_peer_init}
\end{align}
Here,
$p_h(T)=\nabla_y C\big(y_h(T)\big)\T$ and the Jacobians of $F$ are block diagonal matrices $\nabla_YF(Y_n,U_n)=\diag_i\big(\nabla_{Y_{ni}}f(Y_{ni},U_{ni})\big)$ and
$\nabla_UF(Y_n,U_n)=\diag_i\big(\nabla_{U_{ni}}f(Y_{ni},U_{ni})\big)$. The generalized normal cone
mapping $N_{U^s}\big(U_n\big)$ is defined by
\begin{align}
\label{def_cone_peer}
N_{U^s}(u) =&\, \left\{ w\in\R^{sd}: w^T(v-u)\le 0 \mbox{ for all } v\in U_{ad}^s\subset\R^{sd}\right\}.
\end{align}
\par
The discrete KKT conditions \eqref{OCprob_peer_init}-\eqref{KKT_ctr_peer_init} should be good approximations to the continuous ones \eqref{KKT_state}-\eqref{KKT_ctr}.
In what follows, we assume sufficient smoothness of the optimal control problem
such that a local solution $(y^\star,u^\star,p^\star)$ of the KKT system
\eqref{OCprob_objfunc}-\eqref{OCprob_ODEinit} exists. Furthermore, let
the Hamiltonian $H(y,u,p):=p\T f(y,u)$ satisfies a coercivity assumption, which is a strong form of
a second-order condition. Then the first-order optimality conditions are also sufficient \cite{Hager2000}.
If $(y,p)$ is sufficiently close to $(y^\star,p^\star)$,
the control uniqueness property introduced in \cite{Hager2000} yields the existence
of a locally unique minimizer $u=u(y,p)$ of the Hamiltonian over all
$u\in U_{ad}$. Substituting $u$ in terms of $(y,p)$ in \eqref{KKT_state}-\eqref{KKT_adj},
gives then the two-point boundary value problem
\begin{align}
\label{RWPy}
y'(t)=&\,g\big(y(t),p(t)\big),\quad y(0)=y_0,\\
\label{RWPp}
p'(t)=&\,\phi\big(y(t),p(t)\big),\quad p(T)=\nabla_y C\big(y(T)\big)\T,
\end{align}
with the source functions defined by
\begin{align}
g(y,p) :=& f\big(y,u(y,p)\big),
\quad \phi(y,p) := -\nabla_yf\big(y,u(y,p)\big)\T p.
\end{align}
The same arguments apply to the discrete first-order optimality system
\eqref{OCprob_peer_init}-\eqref{KKT_ctr_peer_init}. Substituting
the discrete controls $U_n=U_n(Y_n,P_n)$ in terms of $(Y_n,P_n)$ and defining
\begin{align}
\Phi(Y_n,K_n\T P_n):=&\,\left( \phi(Y_{ni},(K_n\T P_{n})_i)\right)_{i=1}^{s},
\quad G(Y_n,P_n):=\left( g(Y_{ni},P_{ni})\right)_{i=1}^{s},
\end{align}
the approximations for the forward
and adjoint differential equations read in a compact form
\begin{align}
A_0Y_0=&\,a\otimes y_0+b\otimes g\big(y_0,p_h(0)\big)+hK_0G(Y_0,P_0),\label{PMvst}\\[1mm]
A_nY_n=&\,B_nY_{n-1}+hK_nG(Y_n,P_n),\quad 1\le n\le N,\label{PMvnm}\\[1mm]
y_h(T)=&\,(w\T\otimes I)Y_N,\label{Pvlsg}\\[1mm]
p_h(0)=&\,(v\T\otimes I)P_0,\label{PMast}\\[1mm]
A_n\T P_n=&\,B_{n+1}\T P_{n+1}-h\Phi(Y_n,K_n\T P_n),\ 0\le n\le N-1,\label{PManm}\\[1mm]
A_N\T P_N=&\,w\otimes p_h(T)-h\Phi(Y_N,K_N\T P_N),\ n=N.\label{PMaend}
\end{align}
Here, the value of $p_h(0)$ is determined by
an interpolant $p_h(0)=(v\T \otimes I)P_0\approx p(0)$ with $v\in\R^s$ of appropriate order.
In a next step, these discrete equations are now treated as a discretization of the two-point
boundary value problem \eqref{RWPy}-\eqref{RWPp}. We will derive order conditions and give
bounds for the global error.

\section{Error Analysis}
\subsection{Order Conditions}
The local error of the standard Peer method and the starting method is easily analyzed by Taylor expansion of the stage residuals, if the exact ODE solutions are used as stages.
Hence, defining ${\bf y}_n^{(k)}(h\cc):=\big(y^{(k)}(t_n+hc_i)\big)_{i=1}^s,\,k=0,1,$  for the forward Peer method, where $\cc=(c_1,\ldots,c_s)\T$, local order $q_1$ means that
\begin{align}\label{ordvdef}
 A_n{\bf y}_n(h\cc)-B_n{\bf y}_{n-1}(h\cc)-hK_n{\bf y}'(h\cc)=O(h^{q_1}).
\end{align}
In all steps of the Peer triplet, requiring local order $q_1$ for the state variable and order $q_2$ for the adjoint solution leads to the following algebraic conditions from \cite{LangSchmitt2020}.
These conditions depend on the Vandermonde matrix $V_q=\big(\eins,\cc,\cc^2,\ldots,\cc^{q-1})\in\R^{s\times q}$, the Pascal matrix $\PP_q=\Big({j-1\choose i-1}\Big)_{i,j=1}^{q}$ and the nilpotent matrix $\tilde E_q=\big(i\delta_{i+1,j}\big)_{i,j=1}^q$ which commutes with $\PP_q=\exp(\tilde E_q)$.
For the different steps \eqref{PMvst}-\eqref{Pvlsg} and \eqref{PMast}-\eqref{PMaend}, and in the same succession we write down the order conditions from \cite{LangSchmitt2020} when $K_n$ is diagonal.
The forward conditions are
\begin{align}\label{OBvst}
 A_0V_{q_1}=&ae_1\T+be_2\T+K_0V_{q_1}\tilde E_{q_1},\ n=0,
 \\\label{OBvnm}
 A_nV_{q_1}=&B_nV_{q_1}\PP_{q_1}^{-1}+K_nV_{q_1}\tilde E_{q_1},\ 1\le n\le N,
 \\\label{OBvlsg}
 w\T V_{q_1}=&\eins\T,
\end{align}
with the cardinal basis vectors $e_i\in\R^s$, $i=1,\ldots,s$.
The conditions for the adjoint methods are given by
\begin{align}\label{OBast}
 v\T V_{q_2}=&e_1\T,
\\\label{OBanm}
A_n\T V_{q_2}=&B_{n+1}\T V_{q_2}\PP_{q_2}-K_nV_{q_2}\tilde E_{q_2},\ 0\le n\le N-1,
\\\label{OBaend}
A_N\T V_{q_2}=&w\eins\T-K_NV_{q_2}\tilde E_{q_2},\ n=N.
\end{align}

\subsection{Bounds for the Global Error}
In this section, the errors $\check Y_{nj}:=y(t_{nj})-Y_{nj}$, $\check P_{nj}:=p(t_{nj})-P_{nj},$ $n=0,\ldots,N$, $j=1,\ldots,s$, are analyzed.
According to \cite{LangSchmitt2020}, the equation for the errors $\check Y\T=(\check Y_0\T,\ldots,\check Y_N\T)$ and $\check P\T=(\check P_0\T,\ldots,\check P_N\T)$ is a linear system of the form
\begin{align}\label{erreq}
\MM_h\check Z=\tau,
\ \check Z=\begin{pmatrix}\check Y\\\check P\end{pmatrix},
\ \tau=\begin{pmatrix}\tau^Y\\\tau^P\end{pmatrix},
\end{align}
where the matrix $\MM_h$ has a $2\times2$-block structure and $(\tau^Y,\tau^P)$ denote
the corresponding truncation errors.
Deleting all $h$-depending terms from $\MM_h$, the block structure of the remaining matrix $\MM_0$
is given by
\begin{align}\label{MMnull}
 \MM_0=\begin{pmatrix}
  M_{11}\otimes I_m&0\\
  M_{21}\otimes\nabla_{yy}C_N&M_{22}\otimes I_m
 \end{pmatrix}
\end{align}
with $M_{11},M_{21},M_{22}\in\R^{s(N+1)\times s(N+1)}$ and a mean value $\nabla_{yy}C_N\in\R^{m\times m}$
of the symmetric Hessian matrix of $C$.
The index ranges of all three matrices are copied from the numbering of the grid, $0,\ldots,N$, for convenience.
In fact, $M_{21}=(e_N\otimes\eins)(e_N\otimes w)\T$ has rank one only with $e_N=(\delta_{Nj})_{j=0}^N$.
The diagonal blocks of $\MM_0$ are nonsingular and its inverse has the form
\begin{align}\label{MMnullinv}
 \MM_0^{-1}=\begin{pmatrix}
  M_{11}^{-1}\otimes I_m&0\\
  -(M_{22}^{-1}M_{21}M_{11}^{-1})\otimes\nabla_{yy}C_N&M_{22}^{-1}\otimes I_m
 \end{pmatrix}.
\end{align}
The diagonal blocks $M_{11},M_{22}$ have a bi-diagonal block structure with identity matrices $I_s$ in the diagonal.
The individual $s\times s$-blocks of their inverses are easily computed with $M_{11}^{-1}$ having lower triangular block form and $M_{22}^{-1}$ upper triangular block form, with blocks
\begin{align}\label{BDInv}
 (M_{11}^{-1})_{nk}=\bar B_n\cdots\bar B_{k+1},\,k\le n,\quad
 (M_{22}^{-1})_{nk}=\bar B_{n+1}\T\cdots \bar B_k\T,\,k\ge n,
\end{align}
with the abbreviations $\bar B_n:=A_n^{-1}B_n$, $1\le n\le N$ and $\tilde B_{n+1}\T:=(B_{n+1}A_n^{-1})\T$, $0\le n<N$.
Empty products for $k=n$ mean the identity $I_s$.
\par
Defining $\Uu:=h^{-1}(\MM_h-\MM_0)$ and rewriting \eqref{erreq} in fixed-point form
\begin{align}
\check Z=&\,h\MM_0^{-1}\Uu\check Z+\MM_0^{-1}\tau\,,
\end{align}
it has been shown in the proof of Theorem~4.1 of \cite{LangSchmitt2020} for smooth right hand sides $f$ and $h\le h_0$ that
\begin{align}\label{FeAbs}
 \|\check Z\|\le 2\max\{\|M_{11}^{-1}\tau^Y\|,\|M_{22}^{-1}\tau^P\|\}
\end{align}
in suitable norms, where these norms are discussed in more detail in Section~\ref{SSupKv} below.
Moreover, due to the lower triangular block structure of $\MM_0$ the estimate for the error in the state variable may be refined (Lemma~4.2 in \cite{LangSchmitt2020}) to
\begin{align}\label{FeYAb}
 \|\check Y\|\le\|M_{11}^{-1}\tau^Y\|+hL\|\check Z\|
\end{align}
with some constant $L$.
Without additional conditions, estimates of the terms on the right-hand side of \eqref{FeAbs} in the form $\|\check Z\|=O(h^{-1}\|\tau\|)$ lead to the loss of one order in the global error.
However, this loss may be avoided by one additional superconvergence condition on the forward and the adjoint method each, which will be considered next.

\subsection{Superconvergence of the Standard Method}\label{SSupKv}
For $s$-stage Peer methods, global order $s$ may be attained in many cases if other properties of the method have lower priority.
For optimal stiff stability properties like A-stability, however, it may be necessary to sacrifice one order of consistency as in \cite{MontijanoPodhaiskyRandezCalvo2019,SchmittWeiner2017}.
Accordingly, in this paper the order conditions for the standard method are lowered by one compared to the requirements in the recent paper \cite{LangSchmitt2020} to local orders $(s,s-1)$, see Table~\ref{TOBed}.
Still, the higher global orders may be preserved to some extent by the concept of superconvergence which prevents the order reduction in the global error by cancellation of the leading error term.
\par
Superconvergence is essentially based on the observation that the powers of the forward stability matrix $\bar B:=A^{-1}B$ may converge to a rank-one matrix which maps the leading error term of $\tau^Y$ to zero.
This is the case if the eigenvalue 1 of $\bar B$ is isolated.
Indeed, if the eigenvalues $\lambda_j,\,j=1,\ldots,s,$ of the stability matrix $\bar B$ satisfy
\begin{align}\label{EWBed}
 1=\lambda_1>|\lambda_2|\ge\ldots|\lambda_s|,
\end{align}
then its powers $\bar B^n$ converge to the rank-one matrix $\eins\eins\T A$ since $\eins$ and $\eins\T A$ are the right and left eigenvectors having unit inner product $\eins\T A\eins=1$ due to the preconsistency conditions $A^{-1}B\eins=\eins$ and $A\mT B\T \eins=\eins$ of the forward and backward standard Peer method, see \eqref{OBvnm}, \eqref{OBanm}.
The convergence is geometric, i.e.
\begin{align}\label{MKonv}
  \|\bar B^n-\eins\eins\T A\|=\|\big(\bar B-\eins\eins\T A\big)^n\|=O(\gamma^n)\to 0,\,n\to\infty,
\end{align}
for any $\gamma\in(|\lambda_2|,1)$.
Some care has to be taken here since the error estimate \eqref{FeAbs} depends on the existence of special norms satisfying $\|\bar B\|_{X_1}:=\|X_1^{-1}\bar BX_1\|_\infty\!=\!1$, resp. $\|\tilde B\T\|_{X_2}\!=\!1$.
Concentrating on the forward error $\|M_{11}^{-1}\tau^Y\|$,  a first transformation of $\bar B$ is considered with the matrix
\begin{align}\label{XTraf}
 X=\begin{pmatrix}1&-\beta\T\\ \eins_{s-1}&I_{s-1}-\eins_{s-1}\beta\T \end{pmatrix},
 \ X^{-1}=\begin{pmatrix}\beta_1&\beta\T\\-\eins_{s-1}&I_{s-1}
 \end{pmatrix},
\end{align}
where $(\beta_1,\beta\T)=\eins\T A$.
Since $Xe_1=\eins$ and $e_1\T X^{-1}=\eins\T A$, the matrix $X^{-1}\bar BX$ is block-diagonal with the dominating eigenvalue 1 in the first diagonal entry.
Due to \eqref{EWBed} there exists an additional nonsingular matrix $\Xi\in\R^{(s-1)\times(s-1)}$ such that the lower diagonal block of $X^{-1}BX$ has norm smaller one, i.e.
\begin{align}
 \|\Xi^{-1}\big(-\eins_{s-1},I_{s-1}\big)\bar B
 \begin{pmatrix}-\beta\T\\ I_{s-1}-\eins_{s-1}\beta\T \end{pmatrix}\Xi\|_\infty=\gamma<1.
\end{align}
Hence, with the matrix
\begin{align}
X_1:=&\,X\begin{pmatrix}1&\\&\Xi \end{pmatrix}
\end{align}
the required norm is found, satisfying
$\|\bar B\|_{X_1}:=\|X_1^{-1}\bar BX_1\|_\infty=\max\{1,\gamma\}=1$ and $\|\bar B-\eins\eins\T A\|_{X_1}=\gamma<1$ in \eqref{MKonv}.
Using this norm in \eqref{FeAbs} and \eqref{FeYAb} for $\epsilon^Y:=M_{11}^{-1}\tau^Y$, it is seen with \eqref{BDInv} that
\begin{align}\label{SupFe}
 \epsilon^Y_n=&\sum_{k=0}^n \eins(\eins\T A\tau_{n-k}^Y)
 +\underbrace{\sum_{k=0}^n\big(\bar B-\eins\eins\T A\big)^k\tau_{n-k}^Y}_{=O(h^s)},
 \ 0\le n<N. 
\end{align}
Only for $\epsilon_N^Y$ a slight modification is required and the factors in the second sum have to be replaced by $(\bar B_N-\eins\eins\T A)(\bar B-\eins\eins\T A)^{k-1}$ for $k>1$ with norms still of size $O(\gamma^k)$.
Now, in all cases the loss of one order in the first sum in \eqref{SupFe} may be avoided if the leading $O(h^s)$-term of $\tau_{n-k}^Y$ is canceled in the product with the left eigenvector, i.e. if $\eins\T A\tau_{n-k}^Y=O(h^{s+1})$.
An analogous argument may be applied to the second term $\|M_{22}^{-1}\tau^P\|$ in \eqref{FeAbs}.
The adjoint stability matrix $\tilde B\T=(BA^{-1})\T$ possesses the same eigenvalues as $\bar B$ and its leading eigenvectors are also known: $\tilde B\T\eins=\eins$ and $(A\eins)\T\tilde B\T =\eins\T B\T=(A\eins)\T$ by preconsistency and an analogous construction applies.
\par
Under the conditions corresponding to local orders $(s,s-1)$ the leading error terms in $\tau_n^Y=\frac1{s!}h^s\eta_s\otimes y^{(s)}(t_n)+O(h^{s+1})$ and $\tau_n^P=\frac1{(s-1)!}\eta_{s-1}^\ast\otimes p^{(s-1)}(t_n)+O(h^s)$ are given by
\begin{align}\label{etas}
\eta_s=&\,\cc^s-A^{-1}B(\cc-\eins)^s-sA^{-1}K\cc^{s-1},\\[1mm]
\eta_{s-1}^\ast=&\,\cc^{s-1}-A^{-T}B\T(\cc+\eins)^{s-1}+(s-1)A^{-T}K\cc^{s-2}.
\end{align}
\begin{table}[t!]
\begin{center}
\begin{tabular}{|l|c|c|}\hline
 Steps&forward: $q_1=s$&adjoint: $q_2=s-1$\\\hline
 Start, $n=0$&\eqref{OBvst}&\eqref{OBanm},\eqref{AdZB},{$\beta=0$}\\\hline
 Standard, $1\le n<N$&\eqref{OBvnm} &\eqref{OBanm}\\
 Superconvergence  &\eqref{SupKvs} &\eqref{SupKvq}\\
 Compatibility &\eqref{LsbE} &\eqref{LsbE}\\\hline
 Last step&\eqref{OBvnm}, $n=N$&\eqref{OBanm}, $n=N-1$\\
 End point& \eqref{OBvlsg}&\eqref{OBaend},\eqref{AdZB},{$\beta=N$}\\\hline
\end{tabular}
\caption{Combined order conditions for the Peer triplet, including
the compatibility condition \eqref{LsbE} and the condition \eqref{AdZB} for
full matrices $K_N$.}\label{TOBed}
\end{center}
\end{table}
Considering now $(\eins\T A)\eta_s$ in \eqref{SupFe} and similarly  $(A\eins)\T\eta_{s-1}^\ast$ the following result is obtained.
\begin{theorem}\label{TSupKnv}
Let the Peer triplet with $s>1$ stages satisfy the order conditions collected in Table~\ref{TOBed} and let the solutions satisfy $y\in C^{s+1}[0,T]$, $p\in C^s[0,T]$.
Let the coefficients of the standard Peer method satisfy the conditions
\begin{align}\label{SupKvs}
 \eins\T\big(A\cc^s-B(\cc-\eins)^s-sK\cc^{s-1}\big)=&\,0,\\\label{SupKvq}
 \eins\T\big(A\T\cc^{s-1}-B\T(\cc+\eins)^{s-1}+(s-1)K\cc^{s-2}\big)=&\,0,
\end{align}
and let \eqref{EWBed} be satisfied.
Assume, that a Peer solution $(Y\T,P\T)\T$ exists and that $f$ and $C$ have bounded second derivatives.
Then, for stepsizes $h\le h_0$ the error of these solutions is bounded by
\begin{align}\label{Fbeide}
  \|Y_{nj}-y(t_{nj})\|_\infty+ h\|P_{nj}-p(t_{nj})\|_\infty=O(h^{s}),
 \end{align}
$n=0,\ldots,N,\,j=1,\ldots,s$.
\end{theorem}
{\bf Proof:}
Under condition \eqref{SupKvs} the representation \eqref{SupFe} shows that $\|M_{11}^{-1}\tau^Y\|=O(h^s)$.
In the same way follows $\|M_{22}^{-1}\tau^P\|=O(h^{s-1})$ from condition \eqref{SupKvq}.
In \eqref{FeAbs} this leads to a common error $\|\check Z\|=O(h^{s-1})$ which may be refined for the state variable with \eqref{FeYAb} to $\|\check Y\|=O(h^s)$.
\qed
\begin{remark}\label{rem:supconv1}
The estimate \eqref{SupFe} shows that superconvergence may be a fragile property and may be impaired if $|\lambda_2|$ is too close to one, leading to very large values in the bound
\[ \sum_{k=0}^n\|\big(\bar B-\eins\eins\T A\big)^k\|_{X_1}\le\frac{1}{1-\gamma}
\]
for the second term in \eqref{SupFe}.
In fact, numerical tests showed that the value $\gamma\cong|\lambda_2|$ plays a crucial role.
While $|\lambda_2|\doteq0.9$ was appropriate for the 3-stage method \texttt{AP3o32f} which shows order 3 in the tests in Section~\ref{STests}, for two 4-stage methods with $|\lambda_2|\doteq0.9$  superconvergence was not seen in any of our 3 test problems below.
Reducing $|\lambda_2|$ farther with additional searches we found that a value below $\gamma=0.8$ may be safe to achieve order 4 in practice.
Hence, the value $|\lambda_2|$ will be one of the data listed in the properties of the Peer methods developed below.
\end{remark}
\begin{remark}\label{rem:supconv2}
By Theorem~\ref{TSupKnv}, weakening the order conditions from the local order pair $(s+1,s)$ to the present pair $(s,s-1)$ combined with fewer conditions for superconvergence preserves global order $h^s$ for the state variable.
However, it also leads to a more complicated structure of the leading error.
Extending the Taylor expansion for $\tau^Y,\tau^P$ and applying the different bounds, a more detailed representation of the state error may be derived,
\begin{align}\label{FeOrds}\notag
\|\check Y\|\le& h^s\big(\frac{\|\eta_s\|}{(1-\gamma)s!}\|y^{(s)}\|+\frac{|\eins\T A\eta_{s+1}|}{(s+1)!}\|y^{(s+1)}\|\big)\\
 &+ h^s\big(\hat L\|p^{(s)}\|+\hat L\|p^{(s-1)}\|\big) +O(h^{s+1})
\end{align}
with some modified constant $\hat L$.
Obviously, the leading error depends on four different derivatives of the solutions.
Since the dependence on $p$ is rather indirect, we may concentrate on the first line in \eqref{FeOrds}.
Both derivatives there may not be correlated and it may be difficult to choose a reasonable combination of both error constants as objective function in the construction of efficient methods.
Still, in the local error the leading term is obviously $\tau^Y\doteq h^s\frac1{s!}\eta_s\otimes y^{(s)}$,  with $\eta_s$ defined in \eqref{etas}, and it may be propagated through non-linearity and rounding errors.
Hence, in the search for methods,
\begin{align}\label{FeKos}
 err_s:=\frac1{s!}\|\eta_s\|_\infty=\frac1{s!}\|\cc^s-A^{-1}B(\cc-\eins)^s-sA^{-1}K\cc^{s-1}\|_\infty
\end{align}
is used as the leading error constant.
\end{remark}

\subsection{Adjoint Order Conditions for General Matrices $K_n$}\label{SAdO}
The number of order conditions for the boundary methods is so large that they may not be fulfilled with the restriction to diagonal coefficient matrices $K_0,K_N$ for $s\ge3$.
Hence, it is convenient to make the step to full matrices $K_0,K_N$ in the boundary methods and the order conditions for the adjoint schemes have to be derived for this case.
Unfortunately, for such matrices the adjoint schemes \eqref{PMaend} and \eqref{PManm} for $n=0$ have a rather unfamiliar form.
Luckily, the adjoint differential equation $p'=-\nabla_y f(y,u)\T p$ is linear.
We abbreviate the initial value problem for this equation as
\begin{align}\label{JDgl}
 p'(t)=-J(t) p(t),\quad p(T)=p_T,
\end{align}
with $J(t)=\nabla_y f\big(y(t),u(t)\big)\T$
and for some boundary index $\beta\in\{0,N\}$, we consider the matrices $A_\beta=(a_{ij}^{(\beta)})$,
$K_\beta=(\kappa_{ij}^{(\beta)})$.
Starting the analysis with the simpler end step \eqref{PMaend}, we have
\begin{align}\label{JEMeth}
 \sum_{j=1}^s a_{ji}^{(N)}P_{Nj}=&\,w_ip_h(T)+hJ(t_{Ni})\sum_{j=1}^s\kappa_{ji}^{(N)}P_{Nj},\ i=1,\ldots,s,
\end{align}
which is some kind of half-one-leg form since it evaluates the Jacobian $J$ and the solution $p$ at different time points. This step must be analyzed for the linear equation \eqref{JDgl} only.
Expressions for the higher derivatives of the solution $p$ follow easily:
\begin{align}
 p''=\,&(J^2-J')p,\quad p'''=(-J^3+2J'J+JJ'-J'')p.
\end{align}
\begin{lemma}\label{LKvoll}
For a smooth coefficient matrix $J(t)$, the scheme \eqref{JEMeth} for the linear differential equation
\eqref{JDgl} has local order 3 under the conditions
\begin{align}\label{AdOB3}
 A_N\T V_3+K_N\T V_3\tilde E_3=w\eins\T
\end{align}
and with $\beta=N$
\begin{align}\label{AdZB}
 \sum_{i=1\atop i\not=j}^s(c_i-c_j)\kappa_{ij}^{(\beta)}=0,\quad j=1,\ldots,s.
\end{align}
\end{lemma}
{\bf Proof:}
Considering the residual of the scheme with the exact solution $p(t)$, Taylor expansion
at $t_N$ and the Leibniz rule yield
\begin{align*}
 \Delta_i=&\,\sum_{j=1}^s a_{ji}^{(N)}p(t_{Nj})-w_ip(T)-hJ(t_{Ni})\sum_{j=1}^s\kappa_{ji}^{(N)}p(t_{Nj})\\
 =&\,\sum_{k=0}^{q-1}\frac{h^k}{k!}\big(\sum_{j=1}^s a_{ji}^{(N)}c_j^k-w_i\big)p^{(k)}
 -\underbrace{h\sum_{k=0}^{q-2}\frac{h^k}{k!}c_i^kJ^{(k)}
  \sum_{j=1}^s\kappa_{ji}^{(N)}\sum_{\ell=0}^{q-2}\frac{h^\ell}{\ell!}c_j^\ell p^{(\ell)}}_{:=\delta}+O(h^q)
\end{align*}
where all derivatives are evaluated at $t_N$. The second term can be further reformulated as
\begin{align*}
\delta =&\, \sum_{k=1}^{q-1}h^{k}\sum_{\ell=0}^{k-1}
  \sum_{j=1}^s\kappa_{ji}^{(N)}\frac{c_i^\ell c_j^{k-\ell-1}}{\ell!(k-\ell-1)!}J^{(\ell)}p^{(k-\ell-1)}\\
 =&\, \sum_{k=1}^{q-1}h^{k}\sum_{\ell=1}^{k}
  \sum_{j=1}^s\kappa_{ji}^{(N)}\frac{c_i^{\ell-1} c_j^{k-\ell}}{(\ell-1)!(k-\ell)!}J^{(\ell-1)}p^{(k-\ell)}.
\end{align*}
Looking at the factors of $h^0,h^1,h^2$ separately leads to the order conditions
\begin{align*}
h^0:\ 0=&\,A_N\T\eins-w,\\
h^1:\ 0=&\,(\sum_{j=1}^sa_{ji}^{(N)}c_j-w_i) p'-\sum_{j=1}^s\kappa_{ji}^{(N)}Jp
 =(-\sum_{j=1}^sa_{ji}^{(N)}c_j+w_i-\sum_{j=1}^s\kappa_{ji}^{(N)})Jp,\mbox{ i.e.}\\[1mm]
 0=&\,A_N\T\cc+K_N\T\eins-w,\\
h^2:\ 0=&\,\frac12(\sum_ja_{ji}^{(N)}c_j^2-w_i)p''
 -\sum_{j=1}^s\kappa_{ji}^{(N)}\big(c_j Jp'+c_iJ'p\big)\\
 =&\,\frac12(\sum_{j=1}^sc_j^2a_{ji}^{(N)}-w_i+2\sum_{j=1}^sc_j\kappa_{ji}^{(N)})J^2p
 -\frac12(\sum_{j=1}^sc_j^2a_{ji}^{(N)}-w_i+2\sum_{j=1}^s\kappa_{ji}^{(N)}c_i)J'p.
\end{align*}
Cancellation of the factor of $J^2$ requires the condition $0=A_N\T\cc^2+2K_N\T\cc-w$,
 which combines with the $h^0,h^1$-conditions to \eqref{AdOB3}.
The factor of $J'$, however, requires $0=A_N\T\cc^2+2D_cK_N\T\eins -w$ with $D_c=\diag(c_i)$.
Subtracting this expression from the factor of $J^2$ leaves $0=K_N\T\cc-D_cK_N\T\eins$, which corresponds to \eqref{AdZB}.
\qed
\begin{remark} Condition \eqref{AdOB3} is the standard version for diagonal $K_N$ from \cite{LangSchmitt2020}.
Hence, the half-one-leg form of \eqref{JEMeth} introduces $s$ additional conditions \eqref{AdZB}, only, for order 3 while the boundary coefficients $K_\beta$, {$\beta=0,N$,} may now contain $s(s-1)$ additional elements.
In fact, a similar analysis for the step \eqref{PManm} with $n=\beta=0$ reveals again \eqref{AdZB} as the only condition in addition to \eqref{OBanm}.
\end{remark}
%
\section{Existence of Boundary Methods Imposes Restrictions on the Standard Method}\label{SExst}
In the previous paper \cite{LangSchmitt2020}, the combination of forward and adjoint order conditions for the standard method $(A,B,K)$ into one set of equations relating only $A$ and $K$ already gave insight on some background of these methods like the advantages of using symmetric nodes.
It also simplifies the actual construction of methods leading to shorter expressions during the elimination process with algebraic software tools.
For ease of reference, this singular Sylvester type equation is reproduced here,
\begin{align}\label{KmbStdO}
 (V_{q_2}\PP_{q_2})\T A(V_{q_1}\PP_{q_1})-V_{q_2}\T AV_{q_1}
 = (V_{q_2}\PP_{q_2})\T KV_{q_1}\PP_{q_1}\tilde E_{q_1}+(V_{q_2}\tilde E_{q_2})\T KV_{q_1}.
\end{align}
A similar combination of the order conditions for the boundary methods, however, reveals crucial restrictions:
the triplet of methods $(A_0,K_0)$, $(A,B,K)$, $(A_N,B_N,K_N)$ has to be discussed together since boundary methods of appropriate order may not exist for any standard method $(A,B,K)$, only for those satisfying certain compatibility conditions required by the boundary methods.
Knowing these conditions allows to design the standard method alone without the ballast of two more methods with many additional unknowns.
This decoupled construction also greatly reduces the dimension of the search space if methods are optimized by automated search routines.
We start with the discussion of the end method.
\subsection{Combined Conditions for the End Method}
We remind that we now are looking for methods having local order $q_1=s$ and $q_2=s-1$ everywhere, which we abbreviate from now on as $(q_1,q_2)=(s,q)$.
In particular this means that $A\T V_q+K\T V_q\tilde E_q=B\T V_q\PP_q$ for the standard method.
Looking for bottlenecks in the design of these methods, we try to identify crucial necessary conditions and consider the three order conditions for the end method $(A_N,B_N,K_N)$ in combination
\begin{align}\label{RBOrd}
 \begin{array}{rl}
 A_NV_s-B_NV_s\PP_s^{-1}-K_NV_s\tilde E_s&=0,\\[1mm]
 B_N\T V_q\PP_q&=A\T V_q+K\T V_q\tilde E_q=B\T V_q\PP_q,\\[1mm]
 A_N\T V_{q}+K_N\T V_{q}\tilde E_q&=w\eins_q\T.
 \end{array}
\end{align}
From these conditions the matrices $A_N,B_N$ may be eliminated, revealing the first restrictions on $B$.
Here, the singular matrix map
\begin{align}
 \LL_{q,s}:\,\R^{q\times s}\to\R^{q\times s},\quad X\mapsto \tilde E_q\T X+X\tilde E_s
\end{align}
plays a crucial role.
\begin{lemma}\label{LRGlre}
A necessary condition for a boundary method $(A_N,B_N,K_N)$ to satisfy \eqref{RBOrd} is
\begin{align}\label{RGlre}
 \LL_{q,s}\big(V_{q}\T K_N V_s\big) =&\eins_q\eins_s\T-V_q\T BV_s\PP_s^{-1}.
\end{align}
\end{lemma}
{\bf Proof:}
The second condition, $V_q\T B_N=V_q\T B$, in \eqref{RBOrd} leads to the necessary equation
\[  V_q\T A_NV_s-V_q\T K_NV_s\tilde E_s=V_q\T BV_s\PP_s^{-1}\]
due to the first condition.
The transposed third condition
\[ V_{q}\T A_NV_s +\tilde E_q\T V_{q}\T K_NV_s =\eins_q w\T V_s=\eins_q \eins_s\T \]
multiplied by the nonsingular matrix $V_s$ may be used to eliminate $A_N$ and leads to
\[ \tilde E_q\T V_{q}\T K_NV_s+V_q\T K_NV_s\tilde E_s =\eins_q \eins_s\T-V_q\T BV_s\PP_s^{-1}\]
which is the equation \eqref{RGlre} from the statement.\qed
\par
Unfortunately, this lemma leads to several restrictions on the design of the methods due to the properties of the map $\LL_{q,s}$.
Firstly, for diagonal matrices $K_N$ the image $\LL_{q,s}\big(V_{q}\T K_N V_s\big)$ has a very restricted shape.
\begin{lemma}\label{LHankel}
If $K\in\R^{s\times s}$ is a diagonal matrix, then $V_q\T KV_s$ and $\LL_{q,s}\big(V_{q}\T K V_s\big)$ are Hankel matrices with constant entries along anti-diagonals.
\end{lemma}
{\bf Proof:}
With $K=\diag(\kappa_i)$, we have $x_{ij}:=e_i\T(V_q\T KV_s)e_j=\sum_{k=1}^s\kappa_kc_k^{i+j-2}$,
showing Hankel form of $X=(x_{ij})=:(\xi_{i+j-1})$ for $i=1,\ldots,s-1$, $j=1,\ldots,s$.
Now, $e_i\T(\tilde E_q\T X+X\tilde E_s)e_j=(i-1)x_{i-1,j}+x_{i,j-1}(j-1)=(i+j-2)\xi_{i+j-2}$, which shows
again Hankel form.\qed
\par
This lemma means that an end method with diagonal $K_N$ only exists if also $V_q\T BV_s\PP_s^{-1}$ on the right-hand side of \eqref{RGlre} has Hankel structure.
Unfortunately it was observed that for standard methods with definite $K$ this is the case for $q_2\le 2$ only (there exist methods with an explicit stage $\kappa_{33}=0$).
\begin{remark}
Trying to overcome this bottleneck with diagonal matrices $K_N$, one might consider adding additional stages of the end method.
However, using general end nodes $(\hat c_1,\ldots,\hat c_{\hat s})$ with $\hat s\ge s$ does not remove this obstacle.
The corresponding matrix $\LL_{q,s}(\hat V_q\T K_N\hat V_s)$ with appropriate Vandermonde matrices $\hat V$ still has Hankel form.
\end{remark}
But even with a full end matrix $K_N$, Lemma~\ref{LRGlre} and the Fredholm alternative enforce restrictions on the standard method $(A,B,K)$ due to the singularity of $\LL_{q,s}$.
This is discussed for the present situation with $s=4,q=3$, only.
The matrix belonging to the map $\LL_{q,s}$ is $I_s\otimes \tilde E_q\T+\tilde E_s\T\otimes I_q$ and its transpose is $I_s\otimes \tilde E_q+\tilde E_s\otimes I_q$.
Hence, the adjoint of the map $\LL_{q,s}$ is given by
\[ \LL_{q,s}\T:\,\R^{q\times s}\to\R^{q\times s},\quad X\mapsto \tilde E_qX+X\tilde E_s\T.\]
Component-wise the map acts as
\[ \LL_{q,s}\T:\ (x_{ij})\mapsto \big(ix_{i+1,j}+jx_{i,j+1}\big)\]
with elements having indices $i>q$ or $j>s$ being zero.
For $q=3,\,s=4,$ one gets
\[ \LL_{3,4}(X)\T=
 \begin{pmatrix}
 x_{12}+x_{21}& 2x_{13}+x_{22}& 3x_{14}+x_{23}& x_{24}\\
 x_{22}+2x_{31}& 2x_{23}+2x_{32}& 3x_{24}+2x_{33}& 2x_{34}\\
 x_{32}& 2x_{33}& 3x_{34}& 0
 \end{pmatrix}.
\]
It is seen that the kernel of $\LL_{3,4}\T$ has dimension 3 and is given by
\begin{align}\label{KernLT}
 X=\begin{pmatrix}
 \xi_1&\xi_2&\xi_3&0\\
 -\xi_2&-2\xi_3&0&0\\
 \xi_3&0&0&0
 \end{pmatrix}.
\end{align}
In \eqref{RGlre}, the Fredholm condition leads to restrictions on the matrix $B$ from the standard scheme.
However, since the matrix $A$ should have triangular form, it is the more natural variable in the search for good methods and an equivalent reformulation of these conditions for $A$ is of practical interest.
\begin{lemma}\label{LLsbre}
Assume that the standard method $(A,B,K)$ has local order $(s,q)=(4,3)$.
Then, end methods $(A_N,B_N,K_N)$ of order $(s,q)=(4,3)$ only exist if the standard method $(A,B,K)$ satisfies the following set of 3 conditions, either for $B$ or for $A$,
\begin{align}\label{LsbE}
 \left.\begin{array}{r}
 \eins\T B\eins=1\\[1mm]
 \eins\T B\cc-\cc\T B\eins=1\\[1mm]
 \eins\T B\cc^2-2\cc\T B\cc+(\cc^2)\T B\eins=1
 \end{array}\right\}
 \left\{\begin{array}{l}
 1=\eins\T A\eins\\[1mm]
 1=\eins\T A\cc-\cc\T A\eins\\[1mm]
 0=\eins\T A\cc^2-2\cc\T A\cc+(\cc^2)\T A\eins.
 \end{array}\right.
\end{align}
\end{lemma}
{\bf Proof:}
Multiplying equation \eqref{RGlre} by $\PP_s$ from the right and using $\tilde E_s\PP_s=\PP_s\tilde E_s$, an equivalent form is
\[ \LL_{q,s}(V_q\T K_NV_s\PP_s)=\eins_q\eins_s\T \PP_s-V_q\T BV_s=:R,\]
with $\eins_s\T\PP_s=(1,2,4,8,\ldots)$ by the binomial formula.
The Fredholm alternative requires that $tr(X\T R)=0$ for all $X$ from \eqref{KernLT}.
We now frequently use the identities $tr\big((vu\T)R\big)=u\T Rv$ and $Ve_i=\cc^{i-1}$
with $e_i\in\R^s$.
The kernel in \eqref{KernLT} is spanned by 3 basis elements.
The first, $X_1=\bar{e}_1e_1\T$ (with the convention $\bar{e}_i\in\R^q$) leads to
\[ 0\stackrel!=tr(X_1\T R)=\bar{e}_1\T \eins_q\eins_s\T e_1-\bar{e}_1\T V_q\T BV_se_1
=1-\eins_s\T B\eins_s.\]
The second basis element is $X_2=\bar{e}_1e_2\T-\bar{e}_2e_1\T$.
Here $tr(X_2\T\eins_q \eins_s\T\PP_s)=\eins_s\T\PP_se_2-\eins_s\T\PP_se_1=1$ and
\[ tr(X_2\T V_q\T BV_s)= \bar{e}_1\T V_q BV_se_2-\bar{e}_2\T V_q BV_se_1
 =\eins_s\T B\cc-\cc\T B\eins_s.\]
For the third element, $X_3=\bar{e}_1e_3\T -2\bar{e}_2e_2\T+\bar{e}_3e_1\T$, one gets $tr(X_3\T\eins_q \eins_s\T\PP_s)=\eins_s\T\PP_s(e_3-2e_2+e_1)=4-4+1=1$.
The third condition on $B$ is
\[ tr(X_3\T V_q\T BV_s)=\eins_s\T B\cc^2-2\cc\T B\cc+(\cc^2)\T B\eins_s.
\]
The versions for $A$ follow from the order conditions. Let again $\eins:=\eins_s$.
The first columns of \eqref{OBvnm} and \eqref{OBanm} show $B\eins=A\eins$ and $\eins\T B=\eins\T A$, which gives $\eins\T B\eins=\eins\T A\eins=1$.
The second column of \eqref{OBvnm} reads $B\cc=A\cc+A\eins-K\eins$ and leads to $\eins\T A\cc=\eins\T B\cc+\eins\T A\eins-\eins\T K\eins$ showing also $\eins\T K\eins=\eins\T A\eins=1$.
Hence, the second condition in \eqref{LsbE} is equivalent with
\[ 1\stackrel!=\eins\T(B\cc)-\cc\T(B\eins)
 =\eins\T(A\cc+A\eins-K\eins)-\cc\T A\eins
 =\eins\T A\cc-\cc\T A\eins.
\]
In order to show the last equivalence in \eqref{LsbE}, we have to look ahead at the forward condition \eqref{OBvnm} for order 3, which is $B\cc^2=A(\cc^2+2\cc+\eins)-2K(\cc+\eins)$.
This leads to
$\eins\T A\cc^2=\eins\T B\cc^2=\eins\T A(\cc^2+2\cc+\eins)-2\eins\T K(\cc+\eins)$, which is equivalent to $2(\eins\T A \cc-\eins\T K\cc)=1$.
Now, this expression is required in the last reformulation which also uses the second adjoint order condition $B\T\cc=A\T\cc-A\T\eins+K\T \eins$, yielding
\begin{align*}
 1\stackrel!=&\,(\eins\T B)\cc^2+(\cc^2)\T(B\eins)-2\cc\T(B\T\cc)\\
  =&\,\eins\T A\cc^2+(\cc^2)\T A\eins-2\cc\T A\T\cc+\underbrace{2(\cc\T A\T\eins-\cc\T K\T \eins)}_{=1}. \qed
\end{align*}
\begin{remark}\label{RemBA}
It can be shown that only the first condition $\eins\T A\eins=1$ in \eqref{LsbE} is required if the matrix $K$ is diagonal.
This first condition is merely a normalization fixing the free common factor in the class $\{ \alpha\cdot(A,B,K):\,\alpha\not=0\}$ of equivalent methods.
The other two conditions are consequences of the order conditions on the standard method $(A,B,K)$ with diagonal $K$. However, the proof is rather lengthy and very technical and is omitted.
\par
The restrictions \eqref{LsbE} on the standard method also seem to be sufficient with \eqref{AdZB} posing no further restrictions.
In fact, with \eqref{LsbE} the construction of boundary methods was always possible in Section~\ref{SConstr}.
\end{remark}

\subsection{Combined Conditions for the Starting Method}
The starting method has to satisfy only two conditions
\begin{align}\label{OBStrt}
 \begin{array}{rl}
 A_0V_s=&\,ae_1\T+be_2\T+K_0V_s\tilde E_s,\\[1mm]
 A_0\T V_q=&\,B\T V_q\PP_q-K_0\T V_q\tilde E_q.
 \end{array}
\end{align}
The first two columns of the first equation may be solved for $a=A_0V_se_1$ and $b=(A_0V_s-K_0V_s\tilde E_s)e_2$ leading to the reduced conditions
\[ (A_0V_s-K_0V_s\tilde E_s)Q_3=0,\quad Q_3:=I_s-e_1e_1\T-e_2e_2\T.
\]
The presence of the projection $Q_3$ leads to changes in the condition compared to the end method.
\begin{lemma}\label{LRGlli}
A necessary condition for the starting method $(A_0,K_0)$ to satisfy \eqref{OBStrt} is
\begin{align}\label{RGlli}
 \Big(\LL_{q,s}\big(\PP_q^{-T}V_{q}\T K_0 V_s\big)-V_q\T BV_s\Big)Q_3=0.
\end{align}
\end{lemma}
{\bf Proof:}
Transposing the second condition from \eqref{OBStrt} and multiplying with $V_sQ_3$ gives
\[ (V_q\T A_0V_s+\tilde E_q\T V_q\T K_0 V_s-\PP_q\T V_q\T BV_s)Q_3=0,\]
and $V_q\T A_0V_sQ_3$ may now be eliminated from both equations, yielding
\[ (\tilde E_q\T V_q\T K_0 V_s+V_q\T K_0V_s\tilde E_s)Q_3=\PP_q\T V_q\T BV_s Q_3.
\]
Again, $P_q\T$ may be moved to the left side and leads to \eqref{RGlli} since it commutes
with $\tilde E_q\T$. \qed
\par\noindent
The situation is now similar to the one for the end method (also concerning a diagonal form of $K_0$) and we consider again the Fredholm condition.
The matrix belonging to the matrix product $\LL_{q,s}()\cdot Q_3$ is $Q_3\otimes \tilde E_q\T+(\tilde E_sQ_3)\T\otimes I_q$ and it has the transpose $Q_3\otimes \tilde E_q+(\tilde E_sQ_3)\otimes I_q$.
This matrix belongs to the map
\begin{align}\label{LTQ}
 X\mapsto \tilde E_q(XQ_3)+(XQ_3)\tilde E_s\T=\LL_{q,s}\T(XQ_3).
\end{align}
For $q=3,s=4$, images of this map are given by
\[ \LL_{3,4}\T(XQ_3)
 =\begin{pmatrix}
  0& 2x_{13}& x_{23}+3x_{14}& x_{24}\\
  0& 2x_{23}& 2x_{33}+3x_{24}& 2x_{34}\\
  0& 2x_{33}& 3x_{34}& 0
 \end{pmatrix}.
\]
Here, the map $\LL_{3,4}\T$ alone introduces no new kernel elements, the kernel of \eqref{LTQ} coincides with that of the map $X\mapsto XQ_3$ given by matrices of the form $X=\tilde X(I-Q_3),\,\tilde X\in\R^{q\times s}$.
Since the right-hand side of \eqref{RGlli} is $V_q\T BV_sQ_3$, the condition for solvability
\[ tr\big(X\T V_q\T BV_sQ_3\big)=tr\big(\tilde X\T V_q\T BV_sQ_3(I-Q_3)\big)=0
\]is always satisfied since $Q_3(I-Q_3)=0$.
Hence, no additional restrictions on the standard method are introduced by requiring the existence of starting methods.
\section{Construction of Peer Triplets}\label{SConstr}
The construction of Peer triplets requires the solution of the collected order conditions from Table~\ref{TOBed} and additional optimization of stability and error properties.
However, it has been observed that some of these conditions may be related in non-obvious ways, see e.g. Remark~\ref{RemBA}.
This means that the accuracy of numerical solutions may be quite poor due to large and unknown rank deficiencies.
Instead, all order conditions were solved here exactly by algebraic manipulation with rational coefficients as far as possible.
\par
The construction of the triplets was simplified by the compatibility conditions \eqref{LsbE} allowing the isolated construction of the standard method $(A,B,K)$ without the many additional parameters of the boundary methods.
Furthermore, an elimination of the matrix $B$ from forward and adjoint conditions derived in \cite{LangSchmitt2020}, see \eqref{KmbStdO}, reduces the number of parameters in $A,K$ to $s(s+3)/2$ elements with $s-1$ additional parameters from the nodes.
This is so since $(A,B,K)$ is invariant under a common shift of nodes and we chose the increments $d_1=c_2-c_1$, $d_j=c_j-c_2,\,j=3,\ldots,s$ as parameters.
Still, due to the mentioned dependencies between conditions, for $s=4$ a 6-parameter family of methods exists which has been derived explicitly (with quite bulky expressions).
\par
However, optimization of stability properties like $A(\alpha)$-stability or error constants $err_s$ from \eqref{FeKos} was not possible in Maple with 6 free  parameters.
Instead, the algebraic expressions were copied to Matlab scripts for some Monte-Carlo-type search routines.
The resulting coefficients of the standard method were finally approximated by rational expressions and brought back to Maple for the construction of the two boundary methods.
\par
At first glance, having the full 6-parameter family of standard methods at hand may seem to be a good work base.
However, the large dimension of the search space may prevent optimal results with reasonable effort.
This can be seen below, where the restriction to symmetric nodes or singly-implicit methods yielded methods with smaller $err_4$ than automated global searches.
\par
$A(\alpha)$-stability of the method may be checked \cite{LangSchmitt2020} by considering the eigenvalue problem for the stability matrix $(A-zK)^{-1}B$, i.e.
\begin{align}\label{WOK}
Bx=\lambda(A-zK)x\gdw K^{-1}(A-\lambda^{-1}B)x=zx.
\end{align}
as an eigenvalue problem for $z\in\C$ where $|\lambda^{-1}|=1$ runs along the unit circle.
Since we focus on A-stable methods, exact verification of this property would have been preferable, of course, but an algebraic proof of A-stability seemed to be out of reach.
It turned out that the algebraic criterion of the second author \cite{Schmitt2015} is rather restrictive (often corresponding to norm estimates, Lemma~2.8 ibd).
On the other hand, application of the exact Schur criterion \cite{Jackiewicz2009,MontijanoPodhaiskyRandezCalvo2019} is not straight-forward and hardly feasible, since the (rational) coefficients of the stability polynomial are prohibitively large for optimized parameters (dozens of decimal places of the numerators). 
\subsection{Requirements for the Boundary Methods}
Since Lemma~\ref{LLsbre} guarantees the existence of the two boundary methods $(A_0,K_0)$ and $(A_N,B_N,K_N)$, their construction can follow after that of the standard method.
Requirements for these two members of the triplet may also be weakened since they are applied once only.
This relaxation applies to the order conditions as shown in Table~\ref{TOBed}, but also to stability requirements.
Still, the number of conditions at the boundaries is so large that the diagonal respectively triangular forms of $K_0,K_N$ and $A_0,A_N$ have to be sacrificed and replaced by some triangular block structure.
Compared to the computational effort of the complete boundary value problem, the additional complexity of the two boundary steps should not be an issue.
However, for non-diagonal matrices $K_0,K_N$ and $s=4$, the 4 additional one-leg-conditions \eqref{AdZB} have to be obeyed.
\par
Weakened stability requirements mean that the last forward Peer step \eqref{PMvnm} for $n=N$ and the two adjoint Peer steps \eqref{PManm} for $n=N-1$ and $n=0$ need not be A-stable and only nearly zero stable if the corresponding stability matrices have moderate norms.
This argument also applies to the two Runge-Kutta steps without solution output \eqref{PMvst} and \eqref{PMaend}.
However, the implementation of these steps should be numerically safe for stiff problems and arbitrary step sizes.
These steps require the solution of two linear systems with the matrices $A_0-hK_0J_0,\ A_N\T-hJ_N\T K_N\T$ or, rather
\begin{align}\label{EndLGS}
 K_0^{-1}A_0-hJ_0,\quad (K_N^{-1}A_N)\T-hJ_N\T,
\end{align}
where $J_0,J_N$ are block diagonal matrices of Jacobians.
These Jacobians are expected to have absolutely large eigenvalues in the left complex half-plane.
For such eigenvalues, non-singularity of the matrices \eqref{EndLGS} is assured under the following eigenvalue condition:
\begin{align}\label{EndRgl}
 \mu_0:=\min_j\mbox{Re}\lambda_j(K_0^{-1}A_0)>0,\quad
 \mu_N:=\min_j\mbox{Re}\lambda_j(K_N^{-1}A_N)>0.
\end{align}
In Table~\ref{TPTB} the constants $\mu_0,\mu_N$ are displayed for all designed Peer triplets as well as the spectral radii $\varrho(A_N^{-1}B_N)$ and $\varrho(B A_0^{-1}),\,\varrho(B_NA^{-1})$ for the boundary steps of the Peer triplet.
In the search for the boundary methods with their exact algebraic parameterizations, these spectral radii were minimized, and if they were close to one the values $\mu_0,\mu_N>0$ were maximized.

\subsection{$A(\alpha)$-Stable 4-Stage Methods}
Although our focus is on A-stable methods, we also shortly consider $A(\alpha)$-stable methods. We like to consider BDF4 (backward differentiation formulas) as a benchmark, since triplets based on BDF3 were the most efficient ones in \cite{LangSchmitt2020}.
In order to distinguish the different methods, we denote them in the form \texttt{AP}$s$\texttt{o}$q_1q_2$\texttt{aaa}, where \texttt{AP} stands for {\bf A}djoint {\bf P}eer method followed by the stage number and the smallest forward and adjoint {\bf o}rders $q_1$ and $q_2$ in the triplet.
The trailing letters are related to properties of the method like {\bf d}iagonal or {\bf s}ingly diagonal implicitness.
\par
The Peer triplet \texttt{AP4o43bdf} based on BDF4 has equi-spaced nodes $c_i=i/4,\,i=1,\ldots,4$, yielding $w=e_4$. The coefficients of the full triplet are given in Appendix~A.1.
Obviously the method is singly-implicit and its well-known stability angle is $\alpha=73.35^o$.
We also monitor the norm of the zero-stability matrix $\|A^{-1}B\|_\infty$, which may be a measure for the propagation of rounding errors. Its value is $\|A^{-1}B\|_\infty\le5.80$.
Since BDF4 has full global order 4, the error constant from \eqref{FeKos} is $err_4=0$.
Still, the end methods were constructed with the local orders $(4,3)$ according to Table~\ref{TOBed}.
The matrices of the corresponding starting method have a leading $3\times 3$ block and a separated last stage. We abbreviate this as block structure $(3+1)$.
All characteristics of the boundary method are given in Table~\ref{TPTB}.
\par
Motivated by the beneficial properties of Peer methods with symmetric nodes seen in \cite{LangSchmitt2020,SchroederLangWeiner2014}, the nodes of the next triplet with diagonally-implicit standard method were chosen symmetric to a common center, i.e. $c_1+c_4=c_2+c_3$.
Unfortunately, searches for large stability angles with such nodes in the interval $[0,1]$ did not find A-stable methods, but the following method \texttt{AP4o43dif} with {\bf f}lip symmetric nodes and $\alpha=84.00^o$, which is an improvement of 10 degrees over BDF4.
Its coefficients are given in Appendix~A.2.
Although there exist A-stable methods with other nodes in $[0,1]$, this method is of its own interest since its error constant $err_4\doteq2.5$e-3 is surprisingly small.
The node vector of \texttt{AP4o43dif} includes $c_4=1$, leading again to $w=e_4$.
Further properties of the standard method $(A,B,K)$ are $\|A^{-1}B\|_\infty\le 2.01$ and the damping factor $|\lambda_2|=0.26$.
See Table~\ref{TPTB} for the boundary methods.
\subsection{A-stable Methods}
By extensive searches with the full 6-parameter family of diagonally-implicit 4-stage methods many A-stable methods were found even with nodes in $[0,1]$.
In fact, regions with A-stable methods exist in at least 3 of the 8 octants in $(d_1,d_3,d_4)$-space.
Surprisingly, however, for none of these methods the last node $c_4$ was the rightmost one.
Also, it may be unexpected that some of the diagonal elements of $A$ and $K$ have negative signs.
This does not cause stability problems if $a_{ii}\kappa_{ii}>0,\,i=1,\ldots,s$, see also \eqref{EndLGS}.
A-stability assured, the search procedure tried to minimize a linear combination $err_s+\delta\|A^{-1}B\|_\infty$ of the error constant and the norm of the stability matrix with small $\delta<10^{-3}$ to account for the different magnitudes of these data.
As mentioned in Remark~\ref{rem:supconv1}, it was also necessary to include the damping factor $|\lambda_2|$ in the minimization process.
One of the best A-stable standard methods with {\bf g}eneral nodes was named \texttt{AP4o43dig}. Its coefficients are given in Appendix A.3.
It has an error constant $err_4\le 0.0260$ and damping factor $|\lambda_2|\doteq0.798$, see Table~\ref{TPT}.
No block structure was chosen in the boundary methods in order to avoid large norms of the mixed stability matrices $BA_0^{-1},B_NA^{-1},A_N^{-1}B_N$, see Table~\ref{TPTB}.
\par
It was observed that properties of the methods may improve, if the nodes have a wider spread than the standard interval $[0,1]$.
In our setting, the general vector $w$ allows for an end evaluation $y_h(T)$ at any place between the nodes.
Since an evaluation roughly in the middle of all nodes may have good properties, in a further search the nodes were restricted to the interval $[0,2]$.
Indeed, all characteristic data of the method \texttt{AP4o43die} with {\bf e}xtended nodes presented in Appendix~A.4 have improved.
As mentioned before, the standard method is invariant under a common node-shift and a nearly minimal error constant was obtained with the node increments $d_1=\frac{10}{11}$, $d_3=-1$ and $d_4=\frac23$.
Then, the additional freedom in the choice of $c_2$ was needed for the boundary methods, since the conditions \eqref{EndRgl} could only be satisfied in a small interval around $c_2=\frac54$.
The full node vector with this choice has alternating node increments since $c_3<c_1<c_2<c_4$.
The method is A-stable, its error constant $err_4\le 0.0136$ is almost half as large as for the method \texttt{AP4o43dig}, and $\|A^{-1}B\|_\infty\le6.1$ and $|\lambda_2|\le\frac23$ are smaller, too.
The data of the boundary methods can be found in Table~\ref{TPTB}.

\begin{table}
 \begin{tabular}{|c|l|c|c|c|c|c|p{3cm}|}\hline
  $s$&triplet &nodes&$\alpha$&$\|A^{-1}B\|_\infty$&$|\lambda_2|$&$err_s$&Remarks\\\hline
 4& AP4o43bdf&BDF4&$73.35^o$&5.79&0.099&0&singly-implicit\\
  &AP4o43dif& $[0,1]$&$84.0^o$&$2.01$&0.26&0.0025&diag.-implicit\\
  &AP4o43dig&$[0,1]$&$90^o$&24.5&0.798&0.0260&$c_3$=1\\
  &AP4o43sil&$[0,1]$&$90^o$&32.2&0.60& 0.0230&$c_3$=1, sing.impl.\\ 
  &AP4o43die&$[0,2]$&$90^o$&6.08&0.66&0.0135&nodes alternate\\
  \hline
 3&AP3o32f&$[0,1]$&$90^o$&15.3&0.91&0.0170&nodes alternate\\\hline
 \end{tabular}
 \caption{Properties of the standard methods of Peer triplets.}\label{TPT}
\end{table}

\begin{table}
 \begin{tabular}{|l|c|c|c|c|c|c|c|}\hline
 &\multicolumn{3}{c|}{Starting method}&\multicolumn{4}{c|}{End method}\\\cline{2-8}
 triplet &blocks&$\mu_0$&$\varrho(BA_0^{-1})$& blocks&$\mu_N$&$\varrho(A_N^{-1}B_N)$&$\varrho(B_NA^{-1})$\\\hline
 AP4o43bdf& 3+1&5.47&1& 1+3&3.81&1&1.15\\
 AP4o43dif&3+1&6.27&1& 1+3&4.40&1&1.03\\
 AP4o43dig&4&0.99&1& 4&0.89&1.001&1\\
 AP4o43sil&3+1&1.88&1& 4&0.72&1&1.03\\ 
 AP4o43die&3+1&3.80&1& 1+3&0.66&2.6&1.98\\
  \hline
 AP3o32f&1+1+1&1.50&1.02& 1+2&0.94&1&1\\\hline
 \end{tabular}
 \caption{Properties of the boundary methods of Peer triplets.}\label{TPTB}
\end{table}

For medium-sized ODE problems, where direct solvers for the stage equations may be used, an additional helpful feature is diagonal singly-implicitness of the standard method.
In our context this means that the triangular matrices $K^{-1}A$ and $AK^{-1}$ have a constant value $\theta$ in the main diagonal.
Using the ansatz
\[ a_{ii}=\theta\kappa_{ii},\ i=1,\ldots,s,\]
for $A=(a_{ij})$ and $K=(\kappa_{ij})$, the order conditions from Table~\ref{TOBed} 
lead to a cubic equation for $\theta$ with no rational solutions, in general.
In order to avoid pollution of the algebraic elimination through superfluous terms caused by rounding errors, numerical solutions for this cubic equations were not used until the very end.
This means that also the order conditions for the boundary methods were solved with $\theta$ as a free parameter, only the final evaluation of the coefficients in Appendix~A.5 was done with its numerical value.
Also in the Matlab-search for A-stable methods, the cubic equation was solved numerically and it turned out that the {\bf l}argest positive solution gave the best properties.
Hence, this Peer triplet was named \texttt{AP4o43sil}.
For a first candidate with nearly minimal $err_s+\delta\|A^{-1}B\|_\infty$, the damping factor $\gamma\doteq0.89$ was again too close to one to ensure superconvergence in numerical tests, see Remark~\ref{rem:supconv1}.
However, further searches nearby minimizing the damping factor found a better standard method with $|\lambda_2|=0.60$.
Its nodes are $\cc\T=\left(\frac{1}{50},\frac{3}{5},1,\frac{41}{85}\right)$, the diagonal parameter
$\theta=3.34552931287687520$ is the largest zero of the cubic equation
\[112673616\,\theta^3+106686908\,\theta^2-2102637319\,\theta+1621264295=0.\]
Further properties of the standard method are A-stability, nodes in $(0,1]$ with $c_3=1$, norm $\|A^{-1}B\|_\infty=32.2$, and an error constant $err_4=0.0230$.
The end method $(A_N,B_N,K_N)$ has full matrices $A_N,K_N$, see Table~\ref{TPTB}.
\par
For the sake of completeness, we also present an A-stable diagonally-implicit 3-stage method, since in the previous paper \cite{LangSchmitt2020} we could not find such methods with reasonable parameters.
After relaxing the order conditions by using superconvergence, such methods exist.
Applying all conditions for forward order $s=3$ and adjoint order $s-1=2$, there remains a 5-parameter family depending on the node differences $d_1=c_2-c_1$, $d_3=c_3-c_2$ and 3 elements of $A$ or $K$.
A-stable methods exist in all 4 corners of the square $[-\frac12,\frac12]^2$ in the $(d_1,d_3)$-plane, the smallest errors $err_3$ were observed in the second quadrant.
The method \verb+AP3o32f+ with $(d_1,d_3)=(-\frac{5}{27},\frac25)$ has a node vector with $c_3=1$.
The coefficients can be found in Appendix~A.6.
The characteristic data are $err_3\le 0.017$, $\|A^{-1}B\|_\infty\le15.3$, $|\lambda_2|=0.91$.
The starting method is of standard form with lower triangular $A_0$ and diagonal $K_0$.
\par
The main properties of the newly developed Peer triplets are summarized in Table~\ref{TPT} for the standard methods and Table~\ref{TPTB} for the boundary methods.

\section{Numerical Results}\label{STests}
We present numerical results for all Peer triplets listed in Table~\ref{TPT}
and compare them with those obtained for the third-order four-stage one-step W-method
\texttt{Ros3wo} proposed in \cite{LangVerwer2013} which is also A-stable.
All calculations have been done with Matlab-Version R2019a, using the
nonlinear solver \textit{fsolve} to approximate the overall coupled scheme
(\ref{PMvst})--(\ref{PMaend}) with a tolerance $1e\!-\!14$.
To illustrate the rates of convergence, we consider two nonlinear
optimal control problems, the Rayleigh and a controlled motion problem. A
linear wave problem is studied to demonstrate the practical importance of A-stability
for larger time steps.

\subsection{The Rayleigh Problem}
The first problem is taken from \cite{JacobsonMayne1970} and describes the behaviour of
a tunnel-diode oscillator. With the electric current $y_1(t)$ and the
transformed voltage at the generator $u(t)$, the Rayleigh problem reads
\begin{align}
\mbox{Minimize } \int_0^{2.5}\left( u(t)^2+y_1(t)^2\right)\,dt \label{rayleigh_objfunc} &\\
\mbox{subject to } y''_1(t)-y'_1\left( 1.4 - 0.14 y'_1(t)^2 \right)+y_1(t)
=&\,  4u(t),\quad t\in(0,2.5], \label{rayleigh_ODE}\\
y_1(0)=y'_1(0)=&\,-\!5. \label{rayleigh_ODEinit}
\end{align}
Introducing $y_2(t)=y_1'(t)$ and eliminating the control $u(t)$ yields
the following nonlinear boundary value problem (see \cite{LangVerwer2013}
for more details):
\begin{align}
y'_1(t) =&\, y_2(t),\label{rayleigh-bvp1}\\
y'_2(t) =&\, -y_1(t) + y_2\left( 1.4 - 0.14 y_2(t)^2 \right) - 8p_2(t),\\
         &\, y_1(0) = -5,\,y_2(0)=-5,\\
p'_1(t) =&\, p_2(t) - 2 y_1(t),\\
p'_2(t) =&\, -p_1(t) - (1.4-0.42 y_2(t)^2)p_2(t),\\
            &\, p_1(2.5) = 0,\,p_2(2.5)=0.\label{rayleigh-bvp2}
\end{align}
\begin{figure}[t]
\centering
\includegraphics[width=6.8cm]{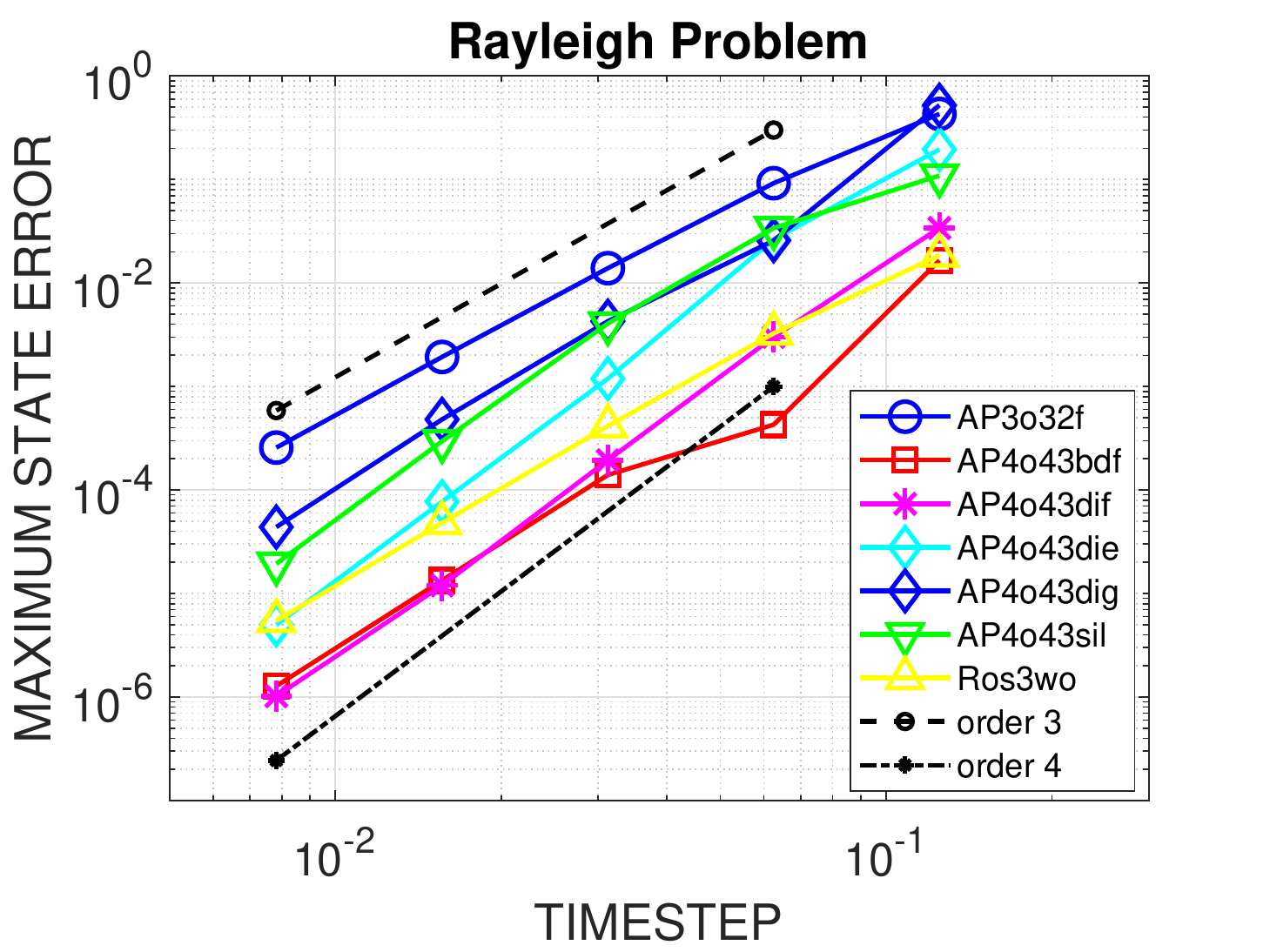}
\hspace{0.1cm}
\includegraphics[width=6.8cm]{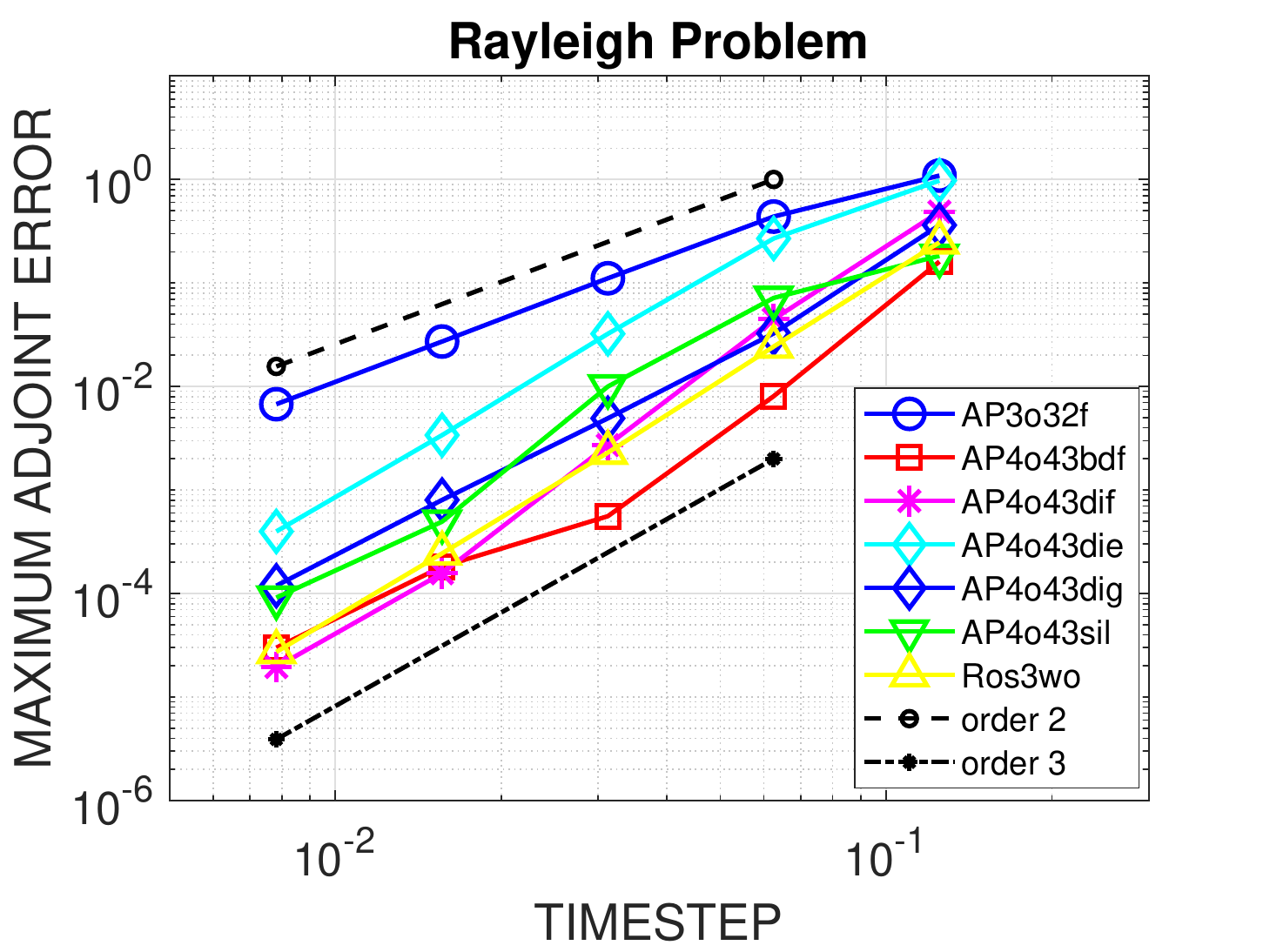}
\parbox{13cm}{
\caption{Rayleigh Problem: Convergence of the maximal state
errors $\|(w\T\otimes I)Y_{n}-y(t_{n+1})\|_\infty$ (left) and
adjoint errors $\|(v\T\otimes I)P_{n}-p(t_n)\|_\infty$ (right), $n=0,\ldots,N$.}
\label{fig:rayleigh}
}
\end{figure}
To study convergence orders of our new methods, we compute a reference solution
at the discrete time points $t=t_n$ by applying the classical fourth-order RK4
with $N\!=\!1280$ steps. Numerical results for the maximum state and adjoint errors are
presented in Figure~\ref{fig:rayleigh} for $N\!+\!1\!=\!20,40,80,160,320$.
\texttt{AP3o32f} and \texttt{Ros3wo} show
their expected orders $(3,2)$ and $(3,3)$ for state and adjoint solutions, respectively.
Order three for the adjoint solutions is achieved by all new four-stage Peer methods.
The smaller error constants of \texttt{AP4o43bdf}, \texttt{AP4o43dif}
and \texttt{Ros3wo} are clearly visible.
The additional superconvergence order four
for the state solutions shows up for \texttt{AP4o43die} and \texttt{AP4o43sil} and nearly for
\texttt{AP4o43dif} and \texttt{AP4o43dig}. \texttt{AP4o43bdf} does not reach its full
order four here, too.

\subsection{A Controlled Motion Problem}
This problem was studied in \cite{LiuFrank2021}. The motion of a damped
oscillator is controlled in a double-well potential, where the control $u(t)$ acts
only on the velocity $y_2(t)$. The optimal control problem reads
\begin{align}
\mbox{Minimize } \frac{\alpha}{2}\| y(6)-y_f\|^2+\frac12\int_0^{6} u(t)^2\,dt \label{ctrmotion_objfunc} &\\
\mbox{subject to } y'_1(t)-y_2(t)=&\,  0,\\
y'_2(t)-y_1(t)+y_1(t)^3+\nu y_2(t)=&\,u(t),\quad t\in(0,6], \label{ctrmotion_ODE}\\
y_1(0)=-1,\,y_2(0)=&\,0, \label{ctrmotion_ODEinit}
\end{align}
where $\nu>0$ is the damping parameter and $y_f$ the target final position. As in
\cite{LiuFrank2021}, we set $\nu=1$, $y_f=(1,0)\T$, and $\alpha=10$.

\begin{figure}[t]
\centering
\includegraphics[width=6.8cm]{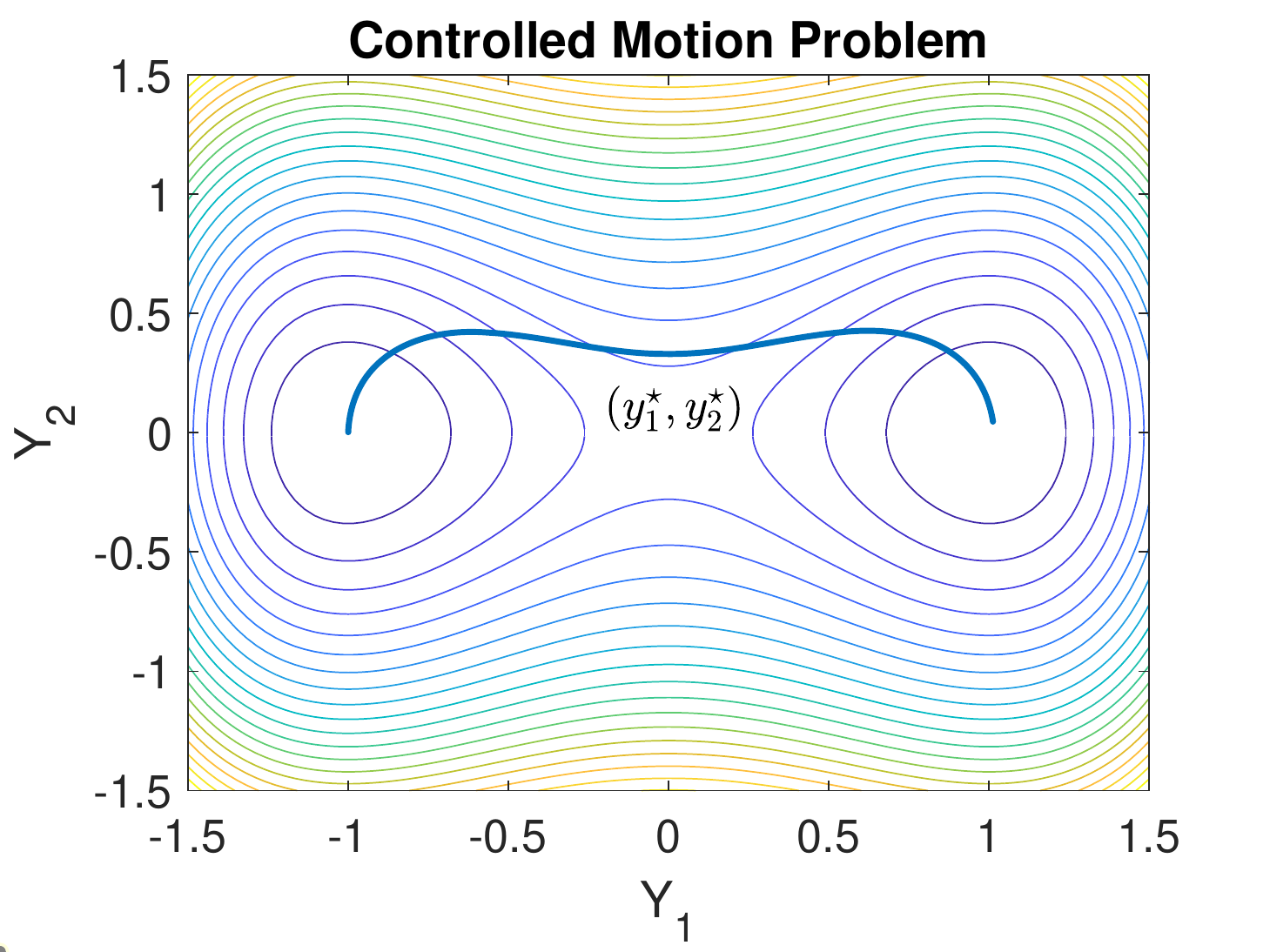}\\
\includegraphics[width=6.8cm]{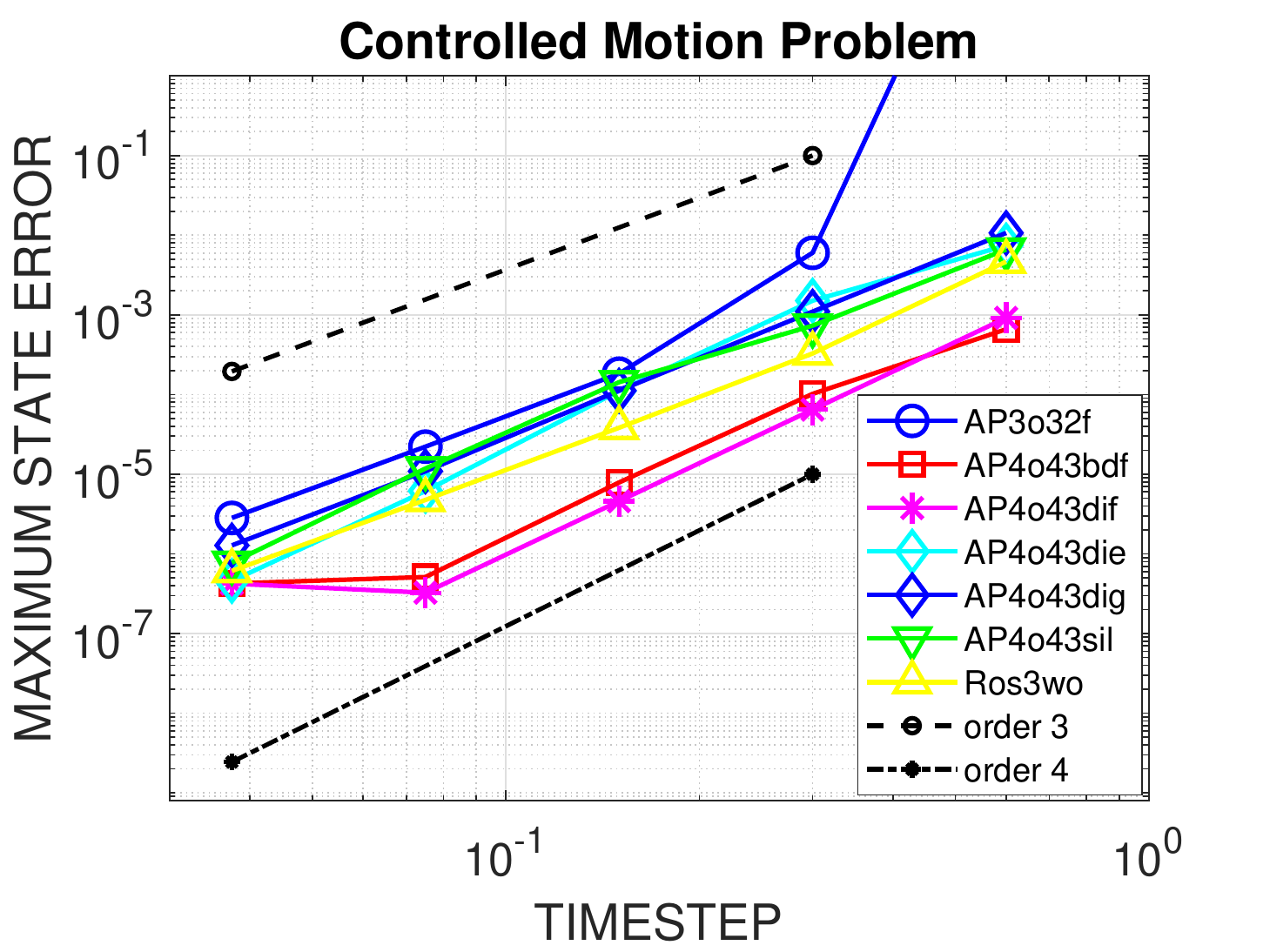}
\hspace{0.1cm}
\includegraphics[width=6.8cm]{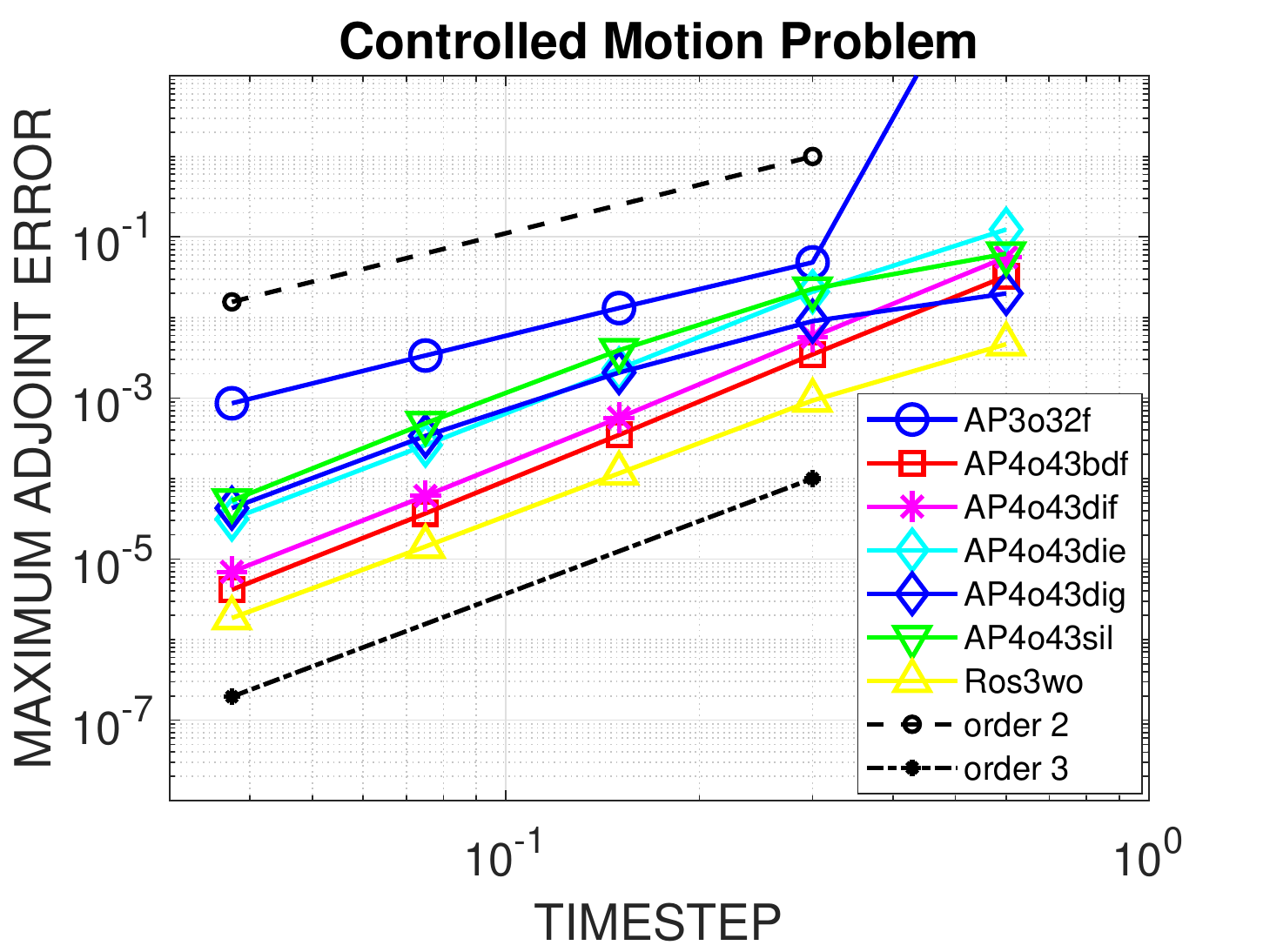}
\parbox{13cm}{
\caption{Controlled Motion Problem: optimal path $(y_1^\star,y_2^\star)$ through the total energy field
$E=\frac12 y_2^2+\frac14y_1^4-\frac12 y_1^2$ visualized by isolines and exhibiting a saddle point
structure (top). Convergence of the maximal state
errors $\|(w\T\otimes I)Y_{n}-y(t_{n+1})\|_\infty$ (bottom left) and
adjoint errors $\|(v\T\otimes I)P_{n}-p(t_n)\|_\infty$ (bottom right),
$n=0,\ldots,N$.}
\label{fig:ctrmotion}
}
\end{figure}

Eliminating the scalar control $u(t)$ yields
the following nonlinear boundary value problem:
\begin{align}
y'_1(t) =&\, y_2(t),\label{ctromotion-bvp1}\\
y'_2(t) =&\, y_1(t) - y_1(t)^3 -\nu y_2(t) - p_2(t),\\
         &\, y_1(0) = -1,\,y_2(0)=0,\\
p'_1(t) =&\, (3y_1(t)^2-1)p_2(t),\\
p'_2(t) =&\, -p_1(t) + \nu p_2(t),\\
            &\, p_1(6) = \alpha (y_1(6)-1),\,p_2(6)=\alpha y_2(6).\label{ctrmotion-bvp2}
\end{align}
The optimal control $u^\star=-p_2^\star$ must accelerate the motion of the particle
to follow an optimal path $(y_1^\star,y_2^\star)$ through the total energy field
$E=\frac12 y_2^2+\frac14y_1^4-\frac12 y_1^2$, shown in Figure~\ref{fig:ctrmotion} on the top,
in order to reach the
final target $y_f$ behind the saddle point. The cost obtained from a reference solution
with $N\!=\!1279$ is $C(y(6))=0.77674$, which
is in good agreement with the lower order approximation in \cite{LiuFrank2021}.
Numerical results for the maximum state and adjoint errors are
presented in Figure~\ref{fig:ctrmotion} for $N\!+\!1\!=\!10,20,40,80,160$.
Worthy of mentioning is the repeated excellent performance of \texttt{AP4o43bdf}
and \texttt{AP4o43dif}, but also the convincing results achieved by the third-order
method \texttt{Ros3wo}. All theoretical orders are well observable, except for \texttt{AP4o43dig},
which tends to order three for the state solutions. A closer inspection reveals that this is
caused by the second state $y_2$, while the first one asymptotically converges with fourth
order. However, the three methods
\texttt{AP4o43die}, \texttt{AP4o43dig}, and \texttt{AP4o43sil} perform quite similar. Observe
that \texttt{AP3o32f} has convergence problems for $N\!=\!9$.

The stagnation of the state errors for the finest step sizes is
due to the limited accuracy of Matlab's \textit{fsolve} -- a fact which was already reported
in \cite{HertyPareschiSteffensen2013}.

\subsection{A Wave Problem}
The third problem is taken from \cite{DontchevHagerVeliov2000} and
demonstrates the practical importance of A-stability. We consider the
optimal control problem
\begin{align}
\mbox{Minimize } y_1(1) + \frac12\,\int_0^1 u(t)^2\,dt \label{wave_objfunc} &\\
\mbox{subject to } y''_1(t)+(2\pi\kappa)^2y_1(t)
=&\,  u(t),\quad t\in(0,1], \label{wave_ODE}\\
y_1(0)=y'_1(0)=&\,0, \label{wave_ODEinit}
\end{align}
where $\kappa=16$ is used. Introducing $y_2(t)=y_1'(t)$ and eliminating the control
$u(t)$ yields the following linear boundary value problem:
\begin{align}
y'_1(t) =&\, y_2(t),\label{wave-bvp1}\\
y'_2(t) =&\, -(2\pi\kappa)^2y_1(t) - p_2(t),\\
         &\, y_1(0) = 0,\,y_2(0)=0,\\
p'_1(t) =&\, (2\pi\kappa)^2p_2(t),\\
p'_2(t) =&\, -p_1(t),\\
            &\, p_1(1) = 1,\,p_2(1)=0.\label{wave-bvp2}
\end{align}
\begin{figure}[t]
\centering
\includegraphics[width=6.8cm]{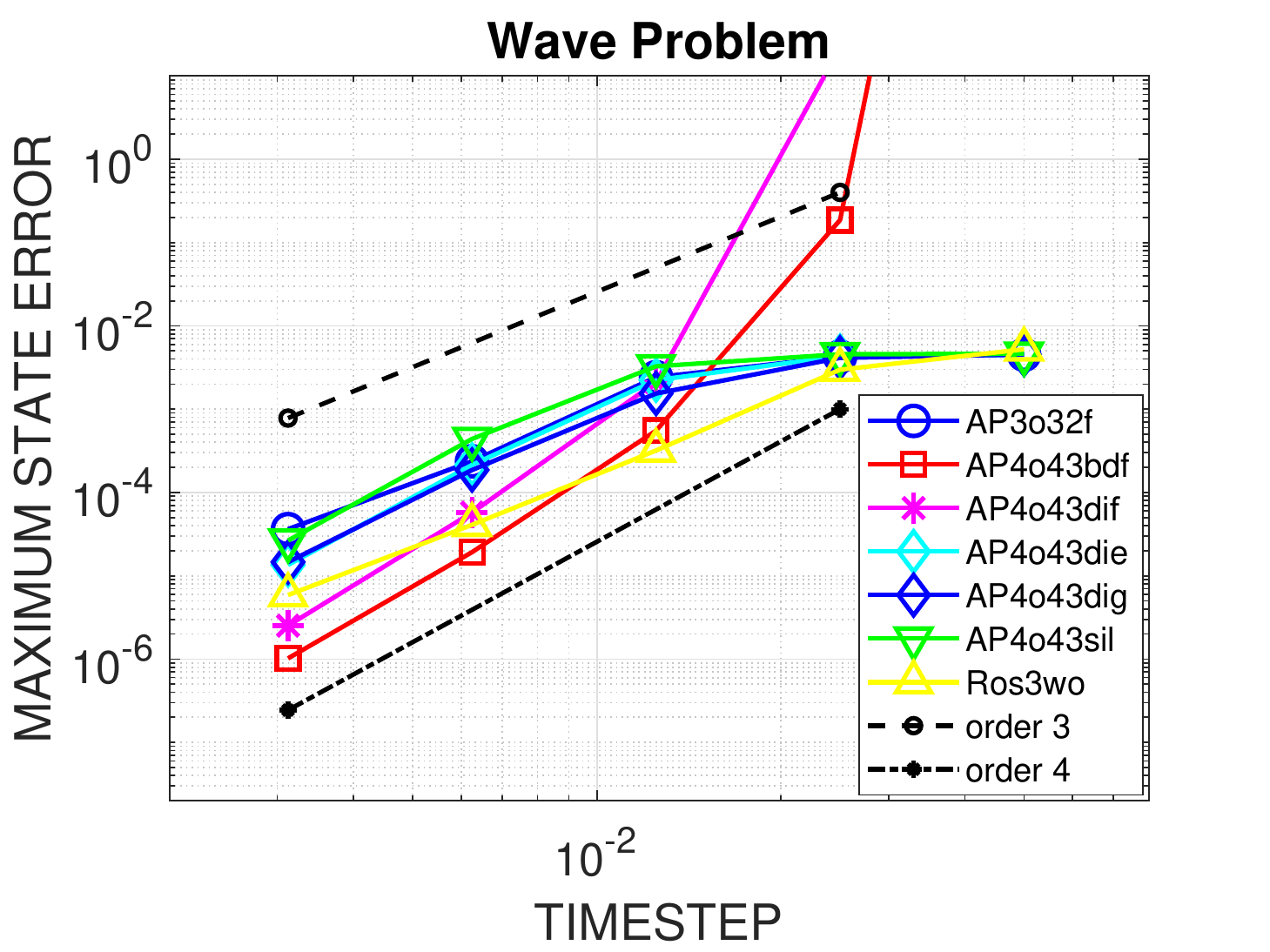}
\hspace{0.1cm}
\includegraphics[width=6.8cm]{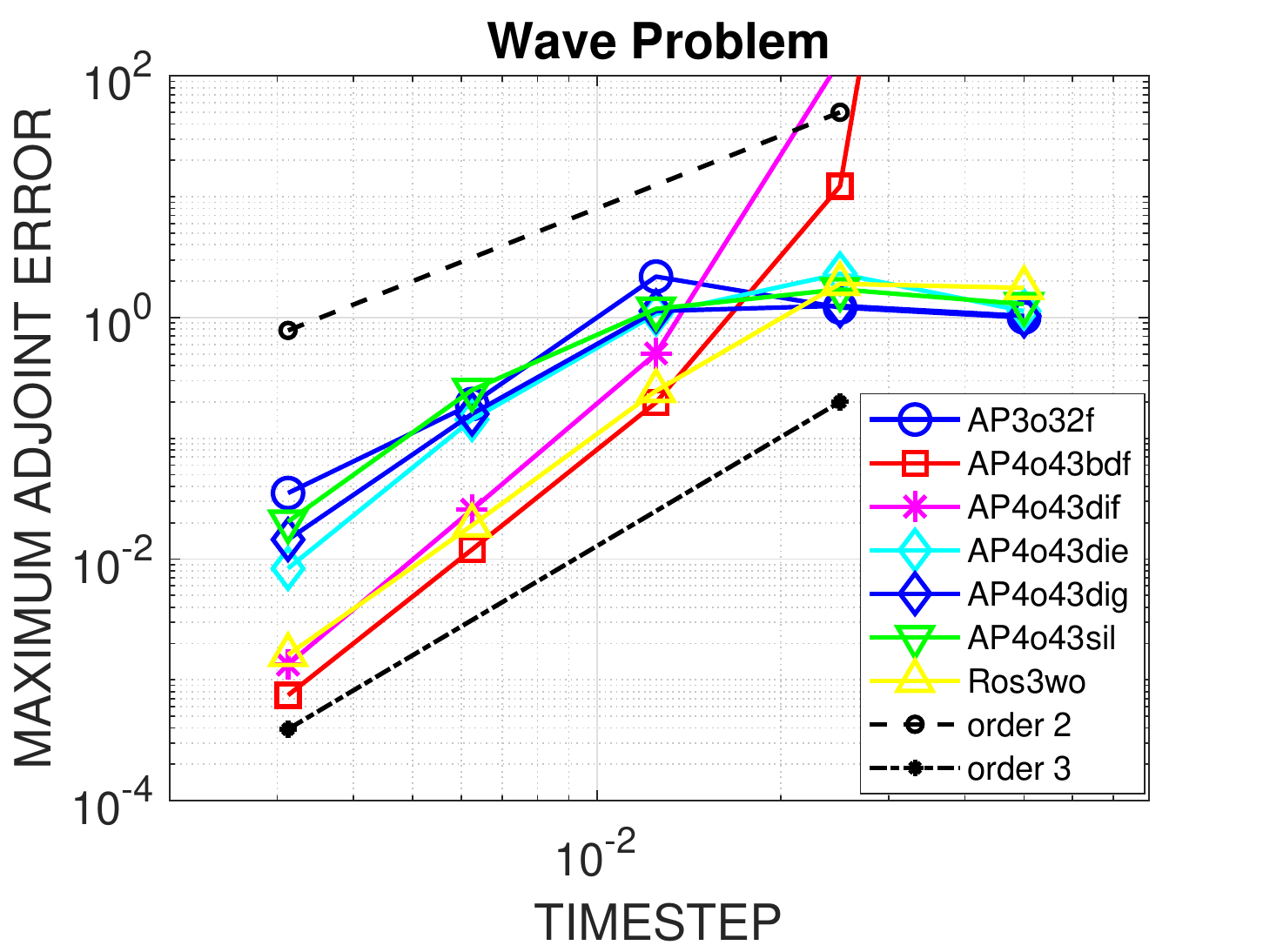}
\parbox{13cm}{
\caption{Wave Problem: Convergence of the maximal state
errors $\|(w\T\otimes I)Y_{n}-y(t_{n+1})\|_\infty$ (left) and
adjoint errors $\|(v\T\otimes I)P_{n}-p(t_n)\|_\infty$ (right),
$n=0,\ldots,N$.}
\label{fig:wave}
}
\end{figure}
The exact solutions are given by
\begin{align}
y_1^*(t)=&\,\frac{1}{2(2\pi\kappa)^3}\sin(2\pi\kappa t)-
\frac{t}{2(2\pi\kappa)^2}\cos(2\pi\kappa t),\\[2mm]
y_2^*(t)=&\,\frac{t}{2(2\pi\kappa)}\sin(2\pi\kappa t),\\[2mm]
p_1^*(t)=&\,\cos(2\pi\kappa t),\;p_2^*(t)=\frac{1}{2\pi\kappa}\sin(2\pi\kappa t),
\end{align}
and the optimal control is $u^*(t)=-p_2^*(t)$. The key observation here is that the
eigenvalues of the Jacobian of the right-hand side in \eqref{wave-bvp1}-\eqref{wave-bvp2}
are $\lambda_{1/2}=2\pi\kappa i$ and $\lambda_{3/4}=-2\pi\kappa i$, which requests
appropriate step sizes for the only $A(\alpha$)-stable methods
\texttt{AP4o43bdf} and \texttt{AP4o34diw} due to their stability restrictions along
the imaginary axis. Indeed, a closer inspection of the stability region of the
(multistep) BDF4 method near the origin reveals that $|\lambda h_{bdf4}|\le 0.3$ is
a minimum requirement to achieve acceptable approximations for problems with imaginary
eigenvalues and moderate time horizon. For the four-stage \texttt{AP4o43bdf}
with step size $h=4h_{bdf4}$, this yields $|\lambda h|\le 1.2$ and hence
$h\le 1/(32\pi)\approx 0.02$
for the wave problem considered. A similar argument applies to \texttt{AP4o34diw}, too.
Numerical results  for the maximum state and adjoint errors are plotted in Figure~\ref{fig:wave} for
$N\!+\!1\!=\!20,40,80,160,320$. They clearly show that both methods deliver first feasible
results for $h=1/80$ and below only, but then again outperform the other Peer methods. Once again,
\texttt{Ros3wo} performs remarkably well. The orders of convergence for the adjoint solutions
are one better than the theoretical values, possibly due to the overall linear structure of
the boundary value problem.

\section{Summary}
We have extended our three-stage adjoint Peer two-step methods constructed in \cite{LangSchmitt2020}
to four stages to not only improve the convergence order of the methods but also their stability.
Combining superconvergence of a standard Peer method with a careful design of starting and end
Peer methods with appropriately enhanced structure, discrete adjoint A-stable Peer methods of order $(4,3)$ could be found.
Still, new requirements had to be dealt with for the higher order pair (4,3).
The property of A-stability comes with larger error constants and a few other, minor
structural disadvantages. As long as A-stability is not an issue to solve the boundary value problem
arising from eliminating the control from the system of KKT conditions, a Peer variant \texttt{AP4o43bdf} of
the $A(73.35^o)$-stable BDF4 is the most attractive method, closely followed by the $A(84^o)$-stable
\texttt{AP4o43dif}, which performs equally well in nearly all of our test cases.
The A-stable methods
\texttt{AP4o43dig} and \texttt{AP4o43die} with diagonally implicit standard Peer methods are very good alternatives if eigenvalues
close or on the imaginary axis are existent. We have also constructed the A-stable method \texttt{AP4o43sil}
with singly-diagonal main part as an additional option if large linear systems can be still
solved by a direct solver and hence the property of requesting one LU-decomposition only is highly valuable.
In future work, we plan to train our novel methods in a projected gradient approach to also tackle
large-scale PDE constrained optimal control problems with semi-discretizations in space.
In these applications, Peer triplets may have to satisfy even more severe requirements.

\vspace{0.5cm}
\par
\noindent {\bf Acknowledgements.}
The first author is supported by the Deutsche Forschungsgemeinschaft
(DFG, German Research Foundation) within the collaborative research center
TRR154 {\em ``Mathematical modeling, simulation and optimisation using
the example of gas networks''} (Project-ID 239904186, TRR154/2-2018, TP B01).

\bibliographystyle{plain}
\bibliography{bibpeeropt}

\begin{thebibliography}{10}

\bibitem{AlbiHertyPareschi2019}
G.~Albi, M.~Herty, and L.~Pareschi.
\newblock Linear multistep methods for optimal control problems and
  applications to hyperbolic relaxation systems.
\newblock {\em Applied Mathematics and Computation}, 354:460--477, 2019.

\bibitem{AlmuslimaniVilmart2021}
I.~Almuslimani and G.~Vilmart.
\newblock Explicit stabilized integrators for stiff optimal control problems.
\newblock {\em SIAM J. Sci. Comput.}, 43:A721--A743, 2021.

\bibitem{BonnansLaurentVarin2006}
F.J. Bonnans and J.~Laurent-Varin.
\newblock Computation of order conditions for symplectic partitioned
  {R}unge-{K}utta schemes with application to optimal control.
\newblock {\em Numer. Math.}, 103:1--10, 2006.

\bibitem{DontchevHagerVeliov2000}
A.L. Dontchev, W.W. Hager, and V.M. Veliov.
\newblock Second-order {R}unge-{K}utta approximations in control constrained
  optimal control.
\newblock {\em SIAM J. Numer. Anal.}, 38:202--226, 2000.

\bibitem{GerischLangPodhaiskyWeiner2009}
A.~Gerisch, J.~Lang, H.~Podhaisky, and R.~Weiner.
\newblock High-order linearly implicit two-step peer – finite element methods
  for time-dependent {PDE}s.
\newblock {\em Appl. Numer. Math.}, 59:624--638, 2009.

\bibitem{GottermeierLang2009}
B.~Gottermeier and J.~Lang.
\newblock Adaptive two-step peer methods for incompressible {N}avier-{S}tokes
  equations.
\newblock In G.~Kreiss, P.~L\"otstedt, A.~Malqvist, and M.~Neytcheva, editors,
  {\em Proceedings of ENUMATH 2009, the 8th European Conference on Numerical
  Mathematics and Advanced Applications, Uppsala, July 2009}, Numerical
  Mathematics and Advanced Applications, pages 387--395. Springer, 2009.

\bibitem{Hager2000}
W.W. Hager.
\newblock {R}unge-{K}utta methods in optimal control and the transformed
  adjoint system.
\newblock {\em Numer. Math.}, 87:247--282, 2000.

\bibitem{HairerWannerLubich2006}
E.~Hairer, G.~Wanner, and Ch. Lubich.
\newblock {\em Geometric Numerical Integration, Structure-Preserving Algorithms
  for Ordinary Differential Equations}, volume~31 of {\em Springer Series in
  Computational Mathematic}.
\newblock Springer, Heidelberg, Berlin, 1970.

\bibitem{HertyPareschiSteffensen2013}
M.~Herty, L.~Pareschi, and S.~Steffensen.
\newblock Implicit-explicit {R}unge-{K}utta schemes for numerical
  discretization of optimal control problems.
\newblock {\em SIAM J. Numer. Anal.}, 51:1875--1899, 2013.

\bibitem{Jackiewicz2009}
Z.~Jackiewicz.
\newblock {\em General Linear Methods for Ordinary Differential Equations}.
\newblock John Wiley\&Sons, Hoboken, New Jersey, 2009.

\bibitem{JacobsonMayne1970}
D.H. Jacobson and D.Q. Mayne.
\newblock {\em Differential Dynamic Programming}.
\newblock American Elsevier Publishing, New York, 1970.

\bibitem{LangSchmitt2020}
J.~Lang and B.A. Schmitt.
\newblock Discrete adjoint implicit peer methods in optimal control.
\newblock Technical Report arXiv:2002.12081, 2020.

\bibitem{LangVerwer2013}
J.~Lang and J.G. Verwer.
\newblock {W}-methods in optimal control.
\newblock {\em Numer. Math.}, 124:337--360, 2013.

\bibitem{LiuFrank2021}
X.~Liu and J.~Frank.
\newblock Symplectic {R}unge-{K}utta discretization of a regularized
  forward-backward sweep iteration for optimal control problems.
\newblock {\em J. Comput. Appl. Math.}, 383:113133, 2021.

\bibitem{LubichOstermann1995}
Ch. Lubich and A.~Ostermann.
\newblock {R}unge-{K}utta approximation of quasi-linear parabolic equations.
\newblock {\em Math. Comp.}, 64:601--627, 1995.

\bibitem{MatsudaMiyatake2021}
T.~Matsuda and Y.~Miyatake.
\newblock Generalization of partitioned {R}unge–{K}utta methods for adjoint
  systems.
\newblock {\em J. Comput. Appl. Math.}, 388:113308, 2021.

\bibitem{MontijanoPodhaiskyRandezCalvo2019}
J.I. Montijano, H.~Podhaisky, L.~Randez, and M.~Calvo.
\newblock A family of {L}-stable singly implicit {P}eer methods for solving
  stiff {IVP}s.
\newblock {\em BIT}, 59:483--502, 2019.

\bibitem{OstermannRoche1992}
A.~Ostermann and M.~Roche.
\newblock {R}unge-{K}utta methods for partial differential equations and
  fractional orders of convergence.
\newblock {\em Math. Comp.}, 59:403--420, 1992.

\bibitem{Sandu2006}
A.~Sandu.
\newblock On the properties of {R}unge-{K}utta discrete adjoints.
\newblock {\em Lecture Notes in Computer Science}, 3394:550--557, 2006.

\bibitem{Sandu2008}
A.~Sandu.
\newblock Reverse automatic differentiation of linear multistep methods.
\newblock In C.~Bischof, H.~B\"ucker, P.~Hovland, U.~Naumann, and J.~Utke,
  editors, {\em Advances in Automatic Differentiation}, volume~64 of {\em
  Lecture Notes in Computational Science and Engineering}, pages 1--12.
  Springer, Berlin, 2008.

\bibitem{SanzSerna2016}
J.M. Sanz-Serna.
\newblock Symplectic {R}unge–{K}utta schemes for adjoint equations, automatic
  differentiation, optimal control, and more.
\newblock {\em SIAM Review}, 58:3--33, 2016.

\bibitem{Schmitt2015}
B.A. Schmitt.
\newblock Algebraic criteria for {A}-stability of peer two-step methods.
\newblock Technical Report arXiv:1506.05738, 2015.

\bibitem{SchmittWeiner2017}
B.A. Schmitt and R.~Weiner.
\newblock Efficient {A}-stable peer two-step methods.
\newblock {\em J. Comput. Appl. Math.}, pages 319--32, 2017.

\bibitem{SchroederLangWeiner2014}
D.~Schr\"oder, J.~Lang, and R.~Weiner.
\newblock Stability and consistency of discrete adjoint implicit peer methods.
\newblock {\em J. Comput. Appl. Math.}, 262:73--86, 2014.

\bibitem{Troutman1996}
J.L. Troutman.
\newblock {\em Variational Calculus and Optimal Control}.
\newblock Springer, New York, 1996.

\end{thebibliography}

\newpage

\section*{Appendix}
In what follows, we will give the coefficient matrices which define the
Peer triplets discussed above. We have used the symbolic option
in Maple as long as possible to avoid any roundoff errors which would pollute
the symbolic manipulations by a great number of superfluous terms. If
possible, we provide exact rational numbers for the coefficients and give
numbers with $16$ digits otherwise. It is sufficient to only show pairs $(A_n,K_n)$ and
the node vector $\cc$, since
all other parameters can be easily computed from the following relations:
\begin{align*}
B_n=&\,(A_nV_s-K_nV_s\tilde E_s)P_sV_s^{-1}, \\[1mm]
a =&\, A_0\eins,\; b=A_0c - K_0\eins, \; w=V_s\mT\eins, \; v=V_s\mT e_1,\quad s=3,4,
\end{align*}
with $e_1=(1,0,\ldots,0)\T\in\R^s$ and the special matrices
\begin{align*}
V_s=\big(\eins,\cc,\cc^2,\ldots,\cc^{s-1}),\quad\PP_q=\Big({j-1\choose i-1}\Big)_{i,j=1}^{s},
\quad\tilde E_s=\big(i\delta_{i+1,j}\big)_{i,j=1}^s\,.
\end{align*}

\subsection*{A.1. Coefficients of \texttt{AP4o54bdf}}\label{coeff:BDF}
\[ \cc\T=\left( \frac14, \frac12, \frac34, 1 \right) \]
\begin{align*}
 A_0=\begin{pmatrix}
  2& \frac12&&\\[2mm]
  -\frac{265}{96}&\frac{17}{96}&\frac{11}{288}&\\[2mm]
  \frac76& -\frac{47}{24}&\frac{25}{12}&\\[2mm]
  -\frac{21}{32}& \frac{227}{96}& -\frac{1163}{288}& \frac{25}{12}
 \end{pmatrix},
 \, K_0=\begin{pmatrix}
  \frac12&&&\\[2mm]
  -\frac{77}{192}& \frac3{32}&&\\[2mm]
   \frac{67}{192}& \frac{17}{96}& \frac{155}{576}&\\[2mm]
   -\frac{19}{192}& -\frac{17}{192}& 0& \frac14
 \end{pmatrix}
\end{align*}
\begin{align*}
 A=\begin{pmatrix}\frac{25}{12}\\[1mm]
  -4&\frac{25}{12}&\\[1mm]
   3&-4&\frac{25}{12}&\\[1mm]
   -\frac43&3&-4&\frac{25}{12}
 \end{pmatrix},
  \, K=\begin{pmatrix}
  \frac14&&&\\[2mm]
  & \frac14&&\\[2mm]
  && \frac14&\\[2mm]
  &&& \frac14
 \end{pmatrix}
\end{align*}
\begin{align*}
 A_N=\begin{pmatrix}
  \frac{635}{96}&&&\\[2mm]
  -\frac{1235}{72}& \frac{35}{32}& \frac{67}{96}& -\frac{43}{288}\\[2mm]
  \frac{4475}{288}& -\frac{35}{24}& 0& \frac{43}{72}\\[2mm]
  -5 & \frac{35}{96}& -\frac{67}{96}& \frac{53}{96}
  \end{pmatrix},
 \, K_N=\begin{pmatrix}
   \frac{25}{32}&&&\\[2mm]
   -\frac53& \frac{61}{192}& -\frac1{192}&\\[2mm]
   \frac{115}{64}& -\frac{13}{48}& \frac{23}{64}&\\[2mm]
   -\frac{185}{288}& \frac{13}{96}& -\frac1{192}& \frac{43}{576}
\end{pmatrix}
\end{align*}

\newpage

\subsection*{A.2. Coefficients of \texttt{AP4o43dif}}
\[ \cc\T=\left( \frac{3}{22}, \frac{53}{132}, \frac{97}{132}, 1\right) \]
\begin{align*}
\small
 A_0=\begin{pmatrix}
  \hfill 1.1582197171362010 & \hfill 0.04624378638947835 & & \\[1mm]
 -0.7020381998219871 & \hfill 1.55936795391810600 & \hfill 0.02624560436867219 & \\[1mm]
 -0.5084723270832399 & -3.71639374160988500 & \hfill 2.21790664582949000 & \\[1mm]
 & \hfill 2.30412222276248700 & -2.66055187655540200 & 1.275350214666080
 \end{pmatrix}
\end{align*}
\begin{align*}
\small
 K_0=\begin{pmatrix}
  \hfill 0.10630513138050800 &&& \\[1mm]
  \hfill 0.39188176473763570 &  \hfill 0.18135555683856680 && \\[1mm]
 -0.05671611789968874 &  \hfill 0.06216151000641002 & 0.3773797141792473 & \\[1mm]
 -0.08101130220826174 & -0.03462160050989925 && 0.1300000000000000
 \end{pmatrix}
\end{align*}
\begin{align*}
\small
 A=\begin{pmatrix}
  \hfill 2.713996187194519 &&& \\[1mm]
 -5.753019558612675 & \hfill 2.063116456071261 && \\[1mm]
  \hfill 5.300000000000000 &-4.392801539381829 & \hfill 2.202673482804081 & \\[1mm]
 -2.313267438350870 & \hfill 2.523025304770754 & -2.619073109161321& 1.275350214666080
 \end{pmatrix}
\end{align*}
\begin{align*}
\small
 K=\diag\left(
 0.2212740342685062, 0.2910929443629617, 0.3576330213685321, 0.13\right)
\end{align*}
\begin{align*}
\small
 A_N=\begin{pmatrix}
  \hfill 3.321208926131899 &&& \\[1mm]
 -6.825690771130220 &  \hfill 1.096545539465253 &  \hfill 0.5537291752348234 & -0.08772005153993665 \\[1mm]
  \hfill 5.504481844998321 & -1.589672639210835 &  \hfill 0.3213889335444567 &  \hfill 0.44690680951897520 \\[1mm]
 -2.000000000000000 &  \hfill 0.493127099745582 & -0.8751181087792801 &  \hfill 0.64081324202096150
  \end{pmatrix}
\end{align*}
\begin{align*}
\small
 K_N=\begin{pmatrix}
  \hfill 0.4780945554021703 &&& \\[1mm]
 -0.7500000000000000 &  \hfill 0.3443968810148793 & 0.03064999258349816 & \\[1mm]
  \hfill 0.9263584006629529 & -0.2474620802680682 & 0.34743387239297130 & \\[1mm]
 -0.4116869618629235 &  \hfill 0.1378269814151266 & 0.03853141924782626 & 0.06599889592052376
\end{pmatrix}
\end{align*}

\newpage

\subsection*{A.3. Coefficients of \texttt{AP4o43dig}}
\[ \cc\T=\left(
\frac{139}{1159}, \frac{11}{19}, 1,
\frac{1375}{2014}\right) \]
\begin{align*}
\small
 A_0=\begin{pmatrix}
 -482.1874750642102 &   4.750000000000000 & -5.916666666666667 &  -6.500000000000000 \\[2mm]
  5295.612100386801 &   78.73229010791468 &  60.16394904407432 &   7.222222222222222 \\[2mm]
  893.8003010294580 &  -4.061724320422766 &  9.736228361340601 &   19.67879886925837 \\[2mm]
 -5707.694317901957 &  -68.66254446706468 &  -58.05689396607474 &  -35.61626763467349
 \end{pmatrix}
\end{align*}
\begin{align*}
\small
 K_0=\begin{pmatrix}
 -49.91295086094522 &  0.5250000000000000 &  -3.439024390243902 &   2.894736842105263 \\[2mm]
  405.5730073881453 &   49.31516975831240 &   8.193548387096774 &  -3.428571428571429 \\[2mm]
  53.62032171809015 &   12.67084977396168 & -1.304859285573478  &  4.013292871986014 \\[2mm]
 -414.6382351541371 & -49.08870633833948 &  -1.334271117642095 &  -15.23643896150272
 \end{pmatrix}
\end{align*}
\begin{align*}
\small
 A=\begin{pmatrix}
 -2.604429828805958 &&& \\[2mm]
  6.603320924494022 &   11.44234275562775 && \\[2mm]
  0.5317173544040980&  -2.710438820206414 &   3.550000000000000 & \\[2mm]
 -5.000000000000000 &   2.026117385005894 &   2.376616772673509 &  -15.21524654319290
 \end{pmatrix}
\end{align*}
\begin{align*}
\small
 K=\diag\left( -0.8973222553064913, 3.337407156628221, 1.164566261468968, -2.604651162790698\right)
\end{align*}
\begin{align*}
\small
 A_N=\begin{pmatrix}
  -3.754385964912281 &   0.01222493887530562 &   1.014925373134328 &  -0.1403508771929825\\[2mm]
  11.35280296428295e &   45.64990373363066   &  -15.17493010383148 &  -20.79910559459065 \\[2mm]
  0.03205698176794144&   2.937595714981687   &  -3.756123160591242 &   2.666624105937791 \\[2mm]
  -7.630473981138614 & -48.59972438748765    &   18.91612789128839 &   18.27283236584584
  \end{pmatrix}
\end{align*}
\begin{align*}
\small
 K_N=\begin{pmatrix}
 -0.9578456075353955 &  -0.3387096774193548 &  0.1194029850746269 &   0.8045112781954887\\[2mm]
  11.57142857142857 &  -96.60667975350700 &  -20.12500000000000 &  102.4846028390543\\[2mm]
  3.761888534906390 &  -33.44512959236514 &  -5.865106921633463 &  34.94704625261853\\[2mm]
 -15.32045224098695 &   134.2026156894527 &   26.37615509001594 & -141.3196101415549
\end{pmatrix}
\end{align*}

\newpage

\subsection*{A.4. Coefficients of \texttt{AP4o43die}}
\[ \cc\T=\left(
\frac{15}{44}, \frac54, \frac14,
\frac{23}{12}\right) \]
\begin{align*}
 A_0=\begin{pmatrix}
\frac{1573}{27} &&& \\[2mm]
-\frac{29140496694667}{1728480384000} & \frac{37071572404007}{69715375488000} & \frac{37246788257}{8450348544} & \\[2mm]
-\frac{98934973036237}{2160600480000} & -\frac{59311823623513}{87144219360000} & -\frac{51770824817}{13203669600} & \\[2mm]
\frac{7653678714559}{1387052160000} & -\frac{7872573544487}{38730764160000} & -\frac{32557703329}{23473190400} & \frac{18014543}{144484600}
 \end{pmatrix}
\end{align*}
\begin{align*}
 K_0=\begin{pmatrix}
\frac{110}{9} &&& \\[2mm]
\frac45 & \frac{168095353644187}{920242956441600} & -\frac{136571614975979}{36809718257664} & \\[2mm]
-\frac{9857504559041}{1080300240000} & \frac{398472962076949}{11503036955520000} & -\frac{18235836500357}{18404859128832} & \\[2mm]
-\frac{1423069729157}{1440400320000} & \frac{398472962076949}{7668691303680000} & \frac{136571614975979}{61349530429440} & \frac{2}{61}
 \end{pmatrix}
\end{align*}
\begin{align*}
 A=\begin{pmatrix}
\frac{45808744223}{19505421000} &&& \\[2mm]
-\frac{279428522}{187552125} & \frac34 && \\[2mm]
\frac{285647}{15004170} & \frac{22704013}{125034750} & -\frac{11}{10} & \\[2mm]
\frac14 & -\frac{11824391}{41678250} & \frac{832579}{4167825} & \frac{18014543}{144484600}
 \end{pmatrix}
\end{align*}
\begin{align*}
 K=\diag\left( \frac{35085281}{25006950}, \frac{2300653}{8335650}, -\frac{1780019}{2500695}, \frac{2}{61}\right)
\end{align*}
\begin{align*}
 A_N=\begin{pmatrix}
\frac{46268635184890481}{6231747944448000} &&& \\[2mm]
-\frac{843924159892681}{239682613248000} & \frac{2095498352743}{2535104563200} & \frac{240331128931}{253510456320} & -\frac{10279031063}{122060590080} \\[2mm]
-4 & \frac{3733202770769}{25351045632000} & -\frac{8883867896017}{6337761408000} & -\frac{192302368031}{24412118016000} \\[2mm]
\frac{39222471881369}{27696657530880} & -\frac{637146509711}{2816782848000} & -\frac{191242568381}{352097856000} & \frac{175009899277}{8137372672000}
  \end{pmatrix}
\end{align*}
\begin{align*}
 K_N=\begin{pmatrix}
\frac{95205971609617}{35952391987200} &&& \\[2mm]
 & \frac{28298016708823}{167316901171200} & -\frac{1296366536717}{4182922529280} & \\[2mm]
-\frac{958020525197}{864240192000} & -\frac{2425439590003}{104573063232000} & -\frac{27877023310129}{41829225292800} & \\[2mm]
-\frac{958020525197}{14980163328000} & -\frac{2425439590003}{69715375488000} & \frac{1296366536717}{6971537548800} & -\frac{22310489177}{4882423603200}
\end{pmatrix}
\end{align*}

\newpage

\subsection*{A.5. Coefficients of \texttt{AP4o43sil}}
\[ \cc\T=\left(\frac{1}{50},\frac35,1,\frac{41}{85} \right) \]
\begin{align*}
\small
 A_0=\begin{pmatrix}
-18.6770976012982273 & -1.15212718448036531 & -0.684527356670693701 & \\
   30.2098963703001422 & -19.0677876392318276 & -7.55433120044842482 & \\
   -9.81986015015644262 & -2.15227598175855777 & 4.86425591259034856 & \\[2mm]
   -\frac{57}{25} & \frac{695}{72} & 8.28572795643498617 & 9.37534909694128499
 \end{pmatrix}
\end{align*}
\begin{align*}
\small
 K_0=\begin{pmatrix}
 -11.4061014637853601 & -0.0776818914116313719 & -0.278650826939386227 &\\[2mm]
   \frac{59}{28} & -13.9738118565057040 & 1.13881886074868390 & \\[2mm]
   -\frac{2498819}{583100} & 2.75133184568842863 & 0.477663277652266390 &\\[2mm]
   \frac{161}{25} & \frac{779}{80} & -0.352459713213735910 & 2.80235150260251923
 \end{pmatrix}
\end{align*}
\begin{align*}
\small
 A=\begin{pmatrix}
   -3.40824065799546119 &&&\\
   -10.5240959029253065 & -6.21116392867196304 &&\\
   1.24215119761892880 & -3.47742608697035235 & 3.58958480260299081 & \\
   12.1231239821473188 & -3.03082301205062472 & 1.32154050930321039 & 9.37534909694123776
 \end{pmatrix}
\end{align*}
\begin{align*}
\small
K=\diag\left( -1.01874482010281778, -1.85655642136068493, 1.07294973886097804, 2.80235150260250919\right)
\end{align*}
\begin{align*}
\small
 A_N=\begin{pmatrix}
 -3.93487127199009570 & \frac{75}{26} & -\frac78 &\\[2mm]
   -4.71932036763822580 & 42.7928650670737006 & -3.60582429888170880 & -36.9240613682294166\\[2mm]
   -\frac{71}{202} & -2.31371287972058488 & 0.191639567731266976 & 1.90723457480523846\\[2mm]
   9.00567678814317298 & -43.3637675719685003 & 5.28918473115044182 & 35.0168267934241781
  \end{pmatrix}
\end{align*}
\begin{align*}
\small
 K_N=\begin{pmatrix}
 -0.687420439097535868 &&&\\[2mm]
   \frac{247}{72} & 24.5972883813771674 & -2.02081177860650990 & -24.7095335501766993\\[2mm]
   -0.427115478996059306 & -\frac{1780}{289} & 0.897376604330415086 & 5.61580307958561348\\[2mm]
   -3.39815952896990200 & -\frac{356}{17} & 1.56153637437775765 & 22.4504633222483055
\end{pmatrix}
\end{align*}

\newpage

\subsection*{A.6. Coefficients of \texttt{AP3o32f}}
\[ \cc\T=\left( \frac{106}{135}, \frac35, 1\right) \]
\begin{align*}
 A_0=\begin{pmatrix}
-\frac{13474483}{2809000} && \\[2mm]
\frac{2765681}{1404500} & \frac{753641}{273375} & \\[2mm]
-\frac{48583191}{81461000} & -\frac{1538339}{1093500} & \frac{1783}{580}
 \end{pmatrix}
\end{align*}
\begin{align*}
 K_0=\diag\left(
 -\frac{13474483}{7155000}, \frac{2513302}{1366875}, \frac{11}{10}\right)
\end{align*}
\begin{align*}
 A=\begin{pmatrix}
 -\frac{11}{2} && \\[2mm]
 \frac{6493}{2700} & \frac{64}{25} & \\[2mm]
 -\frac{25757}{78300} & -\frac{121}{100} & \frac{1783}{580}
 \end{pmatrix}
\end{align*}
\begin{align*}
 K=\diag\left(
 -\frac{93}{50}, \frac{44}{25}, \frac{11}{10}
\right)
\end{align*}
\begin{align*}
 A_N=\begin{pmatrix}
-3 && \\[2mm]
-\frac{559409}{391500} & \frac{5418793}{1458000} & \frac{2257039}{1691280} \\[2mm]
\frac{1733909}{391500} & -\frac{5418793}{1458000} & -\frac{565759}{1691280}
  \end{pmatrix}
\end{align*}
\begin{align*}
 K_N=\diag\left(
  -\frac{1190159}{978750}, \frac{5418793}{3645000}, \frac{2257039}{4228200}
\right)
\end{align*}

\end{document}